\documentclass[twocolumn,9pt]{autart_arxiv} 

\usepackage{times}
\usepackage[T1]{fontenc}
\usepackage[utf8]{inputenc}
\usepackage[english]{babel}
\usepackage{comment}

\usepackage{fancybox,shadow,xcolor} 
\usepackage{setspace}
\usepackage{subfig}
\usepackage{placeins}
\usepackage{interval}
\usepackage{xcolor}
\usepackage{psfrag}
\usepackage{bm}

\usepackage{graphicx}
\usepackage{amssymb,amsmath} %amsthm}
\usepackage{mathrsfs}
\usepackage{mathdots}

\usepackage{enumitem}

\theoremstyle{definition}
\newtheorem{definition}{Definition}

\usepackage{interval}

\DeclareMathOperator*{\argmin}{arg\,min}

\DeclareMathOperator*{\sign}{sign}
\DeclareMathOperator*{\diag}{diag}

\DeclareMathOperator{\Sat}{Sat}

\DeclareMathOperator{\epi}{epi}
\DeclareMathOperator{\dwell}{dwell}
\DeclareMathOperator*{\qd}{\scriptsize qd}

\newcommand{\prox}{\mathbf{prox}}
\newcommand{\dom}{\mathbf{dom}}
\DeclareMathOperator{\Lasso}{L}

\renewcommand{\Re}{\mathbb{R}}

\usepackage{bbm}
\renewcommand{\paragraph}[1]{\smallskip\noindent\textbf{#1.} }

\newcommand{\BM}{\begin{bmatrix}}
\newcommand{\EM}{\end{bmatrix}}

\newcommand{\BBM}{\big(\begin{matrix}}
\newcommand{\EEM}{\end{matrix}\big)}

\newcommand{\bbm}{(\begin{matrix}}
\newcommand{\eem}{\end{matrix})}

\DeclareMathOperator{\dz}{Dz}

\begin{document} 
%\maketitle  

\begin{frontmatter}
																								
\title{Proximal  observers for secure state estimation\thanksref{footnoteinfo}} 

\thanks[footnoteinfo]{This paper was not presented at any IFAC 
meeting. Corresponding author: L. Bako.
}

\author[ECL]{Laurent Bako\thanksref{footnoteinst}}\ead{laurent.bako@centralelille.fr},   
\author[UCBL]{Madiha Nadri}\ead{madiha.nadri-wolf@univ-lyon1.fr}, 
\author[UCBL]{Vincent Andrieu}\ead{vincent.andrieu@gmail.com},              % e-mail address 
\author[INRIA]{Qinghua Zhang}\ead{qinghua.zhang@inria.fr}

\address[ECL]{Univ. Lille, CNRS, Inria, Centrale Lille, UMR 9189 CRIStAL, F-59000 Lille, France}                                              
\address[UCBL]{Universit\'{e} de Lyon, (Universit\'{e} Claude Bernard, CNRS UMR 5007), France}             
\address[INRIA]{Univ. Gustave Eiffel, Inria, Campus de Beaulieu, Rennes, France}  
 
\thanks[footnoteinst]{The work was conducted when this author was with Ecole Centrale de Lyon -- Amp\`{e}re UMR 5005,  France.
}

 \setstretch{.95}         
                                            
\begin{abstract}
This paper discusses a general framework for designing robust state estimators for a class of discrete-time nonlinear systems. We consider  systems that may be impacted by impulsive (sparse but otherwise arbitrary) measurement noise sequences. We show that a family of state estimators, robust to this type of undesired signal, can be obtained by minimizing a class of nonsmooth convex functions at each time step.  
The resulting state observers are defined through proximal operators. We obtain a nonlinear implicit dynamical system in terms of estimation error and prove, in the noise-free setting, that it vanishes asymptotically when  the minimized loss function and the to-be-observed system enjoy appropriate properties.     
From a computational perspective, even though the proposed estimators can be implemented via efficient numerical procedures,  they do not admit closed-form expressions.  The paper argues that by adopting appropriate relaxations, simple and fast analytic expressions can be derived.   
\end{abstract}


\end{frontmatter}

\section{Introduction}
State estimation is a fundamental problem in control theory with decades of developments in both theory \cite{Simon06-Book,Jazwinski70-Book,Bernard22} and applications \cite{Patton97,Mercorelli14,Mercorelli15}. 
In this context, \textit{secure state estimation} refers to the problem of accurately estimating the state of a dynamic system from highly unreliable measurements. We model the unreliability here by an additive signal which is assumed to be sparse but otherwise arbitrary.  
This problem finds application  in the control of cyber-physical systems whose components (sensors, actuators, data processing units, etc.) are linked through some data communication networks. In this context,  the sensor data may be subject to malicious attacks consisting, for example, in the  injection of arbitrary signals to the data traveling over the network. In such cases, the data the state estimator is fed with may be corrupted by arbitrarily large errors (often called outliers).  

Existing relevant works can roughly be divided into three categories. An important body of contributions are batch optimization-based methods which operate on a fixed time horizon over which a finite set of measurements have been collected \cite{Fawzi14-TAC,Farahmand11-TSP,Han19-TAC,Kircher22-TAC}. They mostly apply to discrete-time linear systems. An interesting feature of this estimation scenario is that it can provide estimators having the property of being resilient to the sparse noise \cite{Han19-TAC,Kircher22-TAC}, i.e., with an estimation error which is completely insensitive to the extreme values of the noise.  However, such estimators cannot be directly applied online. 
A second group of methods concern the online estimation of the state.  Some of these use a moving horizon principle to adapt the batch methods for online use purposes; some other types of algorithms run multiple observers in parallel and use a detection rule to isolate bad data \cite{An18-TAC,Mishra16-TNN,Lu19-TAC}. The work presented here is conceptually close to some ideas discussed in \cite{Mattingley10} where a robust Kalman-like estimation strategy is reported. 
Finally, although framed differently in the literature, we note that there is another category of continuous-time methods which are applicable, at least in principle, to secure state estimation  \cite{Alessandri18-Automatica,Astolfi20-TAC,Tarbouriech22-CSL}. The main feature of these latter methods is to pass the measurement through a saturation block before injecting it in the estimator.

The current paper advocates proximal operators as a technical tool which allows for a quite general optimization-based formulation of the online robust state  estimation problem. Note that proximal operators are used in general for solving optimization problems having a composite objective containing a differentiable component and a potentially non smooth one \cite{Parikh14}. Iterative optimization algorithms such as FISTA \cite{Beck09-SIAM} typically rely  on this concept for  optimum searching. In that case, the to-be-optimized objective function is static, that is, independent of time and is fixed over the iterations. This is typically the case in the batch mode of state estimation where  the estimation problem is formed with a fixed and finite collection of data. Here, in contrast, the state estimator is a dynamic system whose dynamics are defined through the solutions of a sequence of time-varying optimization problems. 
Hence, the traditional convergence analysis methods of proximal algorithms do not apply in the setting considered in the current paper. 

\noindent We study a class of adaptive and discrete-time state estimation methods which we call proximal adaptive observers. The principle of such algorithms is to select the state estimate at the current time step by minimizing the distance to the prior  estimate (obtained by applying the predictive model) along with some instantaneous measure relating the  measurement to the model. A general assumption is that the instantaneous cost is convex but not necessarily smooth.  As a consequence, it can be shown, when the output is scalar-valued, that the update equation admits a simple analytic expression involving a  subgradient of the instantaneous cost function. 

\subsection{Outline of the paper}
\noindent We first formulate the general estimation framework in Section \ref{sec:prox-obsv}. 
We then study the stability of the estimators  by considering first the linear case in Section \ref{subsec:stability} and then the nonlinear one in \ref{subsec:analysis-NL}.  
Later in Section \ref{sec:special-cases}, we specialize our framework to some particular choices of the instantaneous cost function, hence yielding different estimators. Those estimators do not admit closed-form expressions in general.  An important contribution of this work is to show that by considering a component-wise relaxation of the ideal measurement update law instead, we can obtain simple and explicit update equations for the estimator. It turns out incidentally that, for three of the candidate loss functions considered in the paper, the expression of the solution involves saturation functions as was suggested  in \cite{Alessandri18-Automatica} in an optimization-free setting. 
The numerical simulation presented in Section \ref{sec:simulation} tends to confirm that our state estimation methods perform very well in the presence of an additive sparse attack signal on the measurements.

\subsection{Notation}
$\Re$ denotes the set of real numbers; $\overline{\Re}=\Re\cup\left\{-\infty,+\infty\right\}$  is the extended real line; $\Re_+$ the set of nonnegative real numbers; $\mathbb{N}$ the set of natural integers and $\mathbb{N}_{\geq n}$ are those natural integers which are greater than or equal to $n$. $\left\|\cdot\right\|_2$ and $|\cdot|$ refer to the Euclidean vector norm and absolute value respectively; $\left\|\cdot\right\|$ is an arbitrary vector norm; $\|\cdot\|_\infty$ is the infinity (vector) norm.  $I$ denotes the square identity matrix of appropriate dimensions. If $W$ is a symmetric positive definite matrix, then $W^{1/2}$ refers to the square root of $W$, i.e., the (unique) symmetric positive definite matrix $S$ such that $S^2=W$. If $S_1$ and $S_2$ are two symmetric matrices, then the notation $S_1\succeq S_2$, or equivalently,  $S_2\preceq S_1$, means that $S_1-S_2$ is positive semidefinite. Similarly $S_1\succ S_2$ (or $S_2\prec S_1$)  signifies that $S_1-S_2$ is positive definite. 

\bigskip

\subsection{Mathematical preliminaries}\label{subsec:maths}
We start by introducing some concepts of optimization theory that will be of key importance in this paper. \\
Consider a function $f:\Re^n\rightarrow \overline{\Re}$ with values on the extended real line. $f$ is called \textit{proper} if its effective domain, $\dom f\triangleq \left\{x\in \Re^n: f(x)<+\infty\right\}$,  is nonempty and $f(x)>-\infty$ for all $x\in \Re^n$; it is called \textit{closed} if its epigraph, $\epi(f)\triangleq \left\{(x,t)\in  \Re^{n}\times \Re: f(x)\leq t\right\}$, is a closed subset of $\Re^{n+1}$. It follows that if $f(x)\in \Re$  for all $x\in \Re^n$, then $f$ is necessarily proper. 
This will essentially be the case for  the to-be-minimized functions considered in the current paper.   \\
 The \textit{proximal operator} of a function $f:\Re^n\rightarrow \overline{\Re}$ is the function $x\mapsto \prox_{f}(x)$ defined on $\Re^n$ by
$$\prox_{f}(x) = \argmin_{z\in\Re^n}\left[\dfrac{1}{2}\left\|z-x\right\|_2^2+f(z)\right]$$
If $f$ is a proper, closed and convex function, then $ \prox_{f}(x)$ exists and is uniquely defined for any $x\in \Re^n$ \cite[Thm 6.3]{Beck17-Book}.   
The \textit{subdifferential} of a  function $f$ at $x\in \dom f$ is the (possibly empty) set $\partial f(x)$ defined by 
$$\partial f(x)=\left\{g: f(y)-f(x)\geq g^\top (y-x) \: \forall y\in  \dom f\right\}. $$
The members of the set $\partial f(x)$ are called \textit{subgradients} of $f$ at $x$. A convex function $f$ is subdifferentiable in the sense that $\partial f(x)\neq \emptyset$ for any  $x$ lying in the interior of $\dom f$ \cite{Beck17-Book}.  Assuming $f$ has a minimum, a point $x^*\in \dom f$ is a minimizer of $f$ if and only if $0\in \partial f(x^*)$.  \\
For more details on the calculations of subdifferentials, we refer, e.g.,  to \cite[Chap. 3]{Beck17-Book}. For our purpose, it is useful to recall the following composition rule for subdifferentials: if $\psi:\Re^n\rightarrow\Re$ is convex, then the function $f$ defined by $f(x)=\alpha\psi(Cx+b)$, for some real matrix $C$, real vector $b$ and nonnegative real $\alpha$, is subdifferentiable and $\partial f(x)=\alpha C^\top \partial\psi(Cx+b)$ for all $x\in \Re^n$.  
%

%%%%%%%%%%%%%%%%%%%%%%%
\section{Secure state estimation problem}

We consider a discrete-time nonlinear and time-varying system described by 
\begin{equation}\label{eq:system}
	\begin{array}{lcl}
		x_{t+1} &=& f_t(x_t)+w_t \\
		 y_t &= &C_t x_t+v_t,
	\end{array}
\end{equation}
where $x_t\in \Re^n$,  $y_t\in \Re^{n_y}$ are the state, the input and the output vectors respectively at time $t\in \mathbb{N}$. As to the variables $w_t$ and $v_t$,  they here represent possible uncertainties of the model. For each time $t$, $C_t\in \Re^{n_y\times n}$ is a matrix and $f_t:\Re^{n}\rightarrow \Re^n$ is a vector field. Note that the system's potential input signals with values  $u_t\in \Re^{n_u}$ are not made explicit in \eqref{eq:system}; this comes without loss of generality since $f_t$ may be taken of the form $f_t(\cdot)={f}_t^\circ(\cdot,u_t)$ for some function ${f}_t^\circ$. 

An informal description of the secure state estimation setting considered in this paper is as follows. The process noise $\left\{w_t\right\}$ is  a dense noise of reasonable amplitude. As to the measurement noise $\left\{v_t\right\}$ however, it can assume arbitrarily large amplitudes. For the sake of an easier description we view the measurement noise $v_t$ as the contribution of two components, 
\begin{equation}\label{eq:decomp-vt}
	v_t = \nu_t+\zeta_t,
\end{equation}
where $\left\{\nu_t\right\}$ is a dense and reasonably bounded noise sequence and  $\left\{\zeta_t\right\}$ is a sparse noise sequence both along the dimension of time and that of the sensors.  $\left\{\zeta_t\right\}$ is unknown (but sparse) and has an impulsive nature, that is, it can take on possibly arbitrarily large values whenever nonzero. It can model, for example, intermittent sensor faults, data losses or adversarial attacks especially when the measurement collection process involves a communication network. 

\paragraph{Problem} 
Given the model $\left\{(C_t,f_t)\right\}$ and the data $\left\{y_t\right\}$ from a model of the form \eqref{eq:system}, we are interested in finding an accurate estimate of the state sequence $\left\{x_t\right\}$ in real-time despite the unknown and possibly arbitrarily large measurement noise $\left\{v_t\right\}$. 
%%%

\section{Proximal robust observers}\label{sec:prox-obsv}

We propose a solution framework to the secure state estimation problem consisting in the  two classical steps of 
 prediction and  measurement update.
Denoting with $\hat{x}_t$ the estimate of the true state $x_t$,   
the prediction step consists in simulating the recurrent equation \eqref{eq:system} based on the state estimate $\hat{x}_{t-1}$ at the previous time step $t-1$, 
\begin{equation}\label{eq:prediction-update}
	\begin{aligned}
	& \hat{x}_{t|t-1} = f_{t-1}(\hat{x}_{t-1}), 
	\end{aligned}
\end{equation}
where the notation $\hat{x}_{t|t-1}$ emphasizes the one step-ahead prediction of $\hat{x}_t$ based on  $\hat{x}_{t-1}$ and the dynamics of the  model. 
The  measurement update step  integrates  the information contained in the measurement $y_t$ into the prior estimate $\hat{x}_{t|t-1}$ through solving an optimization problem as follows  
\begin{equation}\label{eq:measurement-update}
	\begin{aligned}
	& \hat{x}_{t}=\argmin_{z\in \Re^n}\left[\dfrac{1}{2}\big\|W_t^{-1}\left(z-	\hat{x}_{t|t-1}\right)\big\|_2^2+g_t(z)\right]. 
	\end{aligned}
\end{equation}
Here, $W_t\in \Re^{n\times n}$ is a user-defined weighting matrix which we assume to be symmetric and positive definite and  $g_t:\Re^n\rightarrow \Re_+$ is a time-varying function such that for all $z\in \Re^n$, $g_t(z)$ has the form 
\begin{equation}\label{eq:gt}
	g_t(z) = \psi\big(V_t^{-1}(y_t-C_t z)\big)
\end{equation}
 for some convex function $\psi:\Re^{n_y}\rightarrow\Re_+$ and some symmetric and positive definite weighting matrix $V_t\in \Re^{n_y\times n_y}$; $C_t$ is the observation matrix from \eqref{eq:system}.  
The rationale behind the update law \eqref{eq:measurement-update} is to select the estimate $\hat{x}_t$ to be consistent with both the prediction and the measurement models. The weighting matrices $W_t$ and $V_t$ are intended for somewhat normalizing the vector-valued error involved in the to-be-minimized cost function.

Note that the definition \eqref{eq:measurement-update} is indeed equivalent to 
\begin{equation}\label{eq:weigthed-proxi}
	\hat{x}_{t}=W_t\: \prox_{g_t^w}\big(W_t^{-1}\hat{x}_{t|t-1}\big), 
\end{equation}
with $g_t^w$ being the function given by 
$$g_t^w(z)=g_t(W_t z)=\psi\big(V_t^{-1}(y_t-C_t W_tz)\big).$$ 
In case $W_t$ is taken equal to the identity matrix, we get $\hat{x}_{t}=\prox_{g_t}\big(\hat{x}_{t|t-1}\big)$.  
While these two-steps are  quite classical in optimal observer design, the contribution of this work resides more in the question of what would be a good choice for the loss function $\psi$ mentioned above so that the obtained observer is robust to  the sparse/impulsive noise $\zeta_t$.

With the  setting described above,  it follows from the preliminary section that the estimate $\hat{x}_t$ in \eqref{eq:measurement-update} is well-defined in the sense that the minimizer always exists and  is unique. However, no closed form expression of the estimate can be hoped for in general. But since the optimization problem in \eqref{eq:measurement-update} is  convex, it can be numerically solved by efficient tools. As will be shown in Section \ref{sec:special-cases},  with an appropriate implementation strategy, we can even derive fast closed-form solutions. 

For now, let us express the dynamics of the estimation error $e_t\triangleq \hat{x}_t-x_t$ associated with \eqref{eq:prediction-update}-\eqref{eq:measurement-update}. For this purpose, note first that the estimate can be expressed as
$$\hat{x}_t = f_{t-1}(\hat{x}_{t-1})-W_t^2\xi_t(e_t),  $$
with $\xi_t(e_t) \in \partial g_t(\hat{x}_t)=\partial g_t({x}_t+e_t)$.  
Comparing with \eqref{eq:system} gives the error dynamics in the form
\begin{equation}\label{eq:error}
	\begin{aligned}
		e_t=&f_{t-1}({x}_{t-1}+e_{t-1})-f_{t-1}(x_{t-1})\\
		&\qquad \qquad \qquad \qquad  -w_{t-1}-W_t^2\xi_t(e_t). 
	\end{aligned}
\end{equation}

\subsection{Stability analysis: linear case}\label{subsec:stability}
In this section, we analyze  the estimation error \eqref{eq:error} induced by the estimator \eqref{eq:prediction-update}-\eqref{eq:measurement-update}. For clarity reasons, we first study a linear time-varying setting. The nonlinear scenario is deferred to Section \ref{subsec:analysis-NL}.  
Furthermore, to provide a rigorous foundation for the proposed framework, the theoretical analysis is deliberately structured in stages. We first establish the asymptotic convergence of the estimation error in a nominal, noise-free setting ($w_t=0$, $v_t=0$). Proving this nominal stability is a fundamental prerequisite to ensure the structural consistency of the proposed proximal observer, confirming that it accurately recovers the state under ideal conditions. Once this nominal convergence is guaranteed, we subsequently address the robust performance of the estimator in the presence of uncertainties, demonstrating that the estimation error remains bounded even when the measurements are corrupted by arbitrarily large impulsive noise.
\\
Considering systems of the form \eqref{eq:system} where the functions $\left\{f_t\right\}$  are expressible in the form
\begin{equation}\label{eq:linear-ft}
	f_t(x)=A_tx+B_tu_t,
\end{equation}
 with $A_t$ and $B_t$ being matrices of appropriate dimensions, the dynamics of the estimation error \eqref{eq:error} simplifies to
\begin{equation}\label{eq:error-dynamics-linear}
	e_t = A_{t-1}e_{t-1}-w_{t-1} -W_t^2\xi_t(e_t).
\end{equation}
We can then obtain a general convergence result of the estimation error towards zero in a noise-free setting. 
Note that an important challenge of the analysis is that the error dynamics in \eqref{eq:error-dynamics-linear} is an implicit system due to the presence of $e_t$ on the right hand side. 
To present the stability  result,  let us  introduce the following notations for $e\in\Re^n$:  
\begin{align}
&	D_t(e)= e^\top \left(A_{t-1}^\top W_{t}^{-2}A_{t-1}-W_{t-1}^{-2}\right)e \notag\\
                  &  \hspace{4.5cm}-2\psi(V_{t-1}^{-1}C_{t-1}e), \label{eq:Dt}\\
& G_t(e)=e^\top W_t^{-2}e+2\psi(V_t^{-1}C_te)+\xi_t(e)^\top W_t^2\xi_t(e), \label{eq:Gt}\\
& \Sigma_t(e) =\sum_{k=t-T}^{t-1}\xi_k(e_k)^\top W_k^2\xi_k(e_k) \notag \\
             & \qquad \qquad  \mbox{subject to }  e_k=A_{k-1}e_{k-1}-W_k^2\xi_k(e_k), \label{eq:Sigmat}\\
						 & \qquad \qquad \qquad \qquad \: \:   e_{t-T}=e. \notag
\end{align}
where $T\in \mathbb{N}$ is a fixed time horizon and $\xi_t(e)\in \partial \tilde{g}_t(e)$ 
with $\tilde{g}_t(e)=\psi(V_t^{-1}C_te)$ and hence $\partial \tilde{g}_t(e)=C_t^\top V_t^{-1} \partial\psi(V_t^{-1}C_te)$. 
The expressions \eqref{eq:Dt} and  \eqref{eq:Gt} are valid for $t\geq 1$ while \eqref{eq:Sigmat} is for $t\geq T$. 
We also define $G_0$ by $G_0(e)=e^\top W_0^{-2}e$.

We note that the dynamics in \eqref{eq:error-dynamics-linear} and \eqref{eq:Sigmat} are well-defined in the sense that $\xi_k(e_k)$ is uniquely determined even though $\partial \tilde{g}_k(e_k)$ may represent a set.  This is a direct consequence of the lemma below where $\xi_k(e_k)$ satisfies 
$$\xi_k(e_k)\in \partial \tilde{g}_k\big(A_{k-1}e_{k-1}-W_k^2\xi_k(e_k)\big).$$
%%% 
\begin{lem}\label{lem:uniqueness}
Let $f:\Re^n\rightarrow\Re$ be a proper, closed, convex function with $\dom f=\Re^n$ and let $W\in \Re^{n\times n}$ be a symmetric positive definite matrix. Then
\begin{enumerate}
	\item[(a)] For all $x\in\Re^n$, there is a unique $\xi\in \Re^n$ satisfying $\xi\in \partial f(x-W\xi)$. 
Moreover, $\xi=0$ is a solution to this equation if and only if $x\in \argmin_{z}f(z)$. \label{lem-(a)}
\item[(b)] The function $x\mapsto \xi(x)$ such that $\xi(x)\in \partial f\big(x-W\xi(x)\big)$ is well-defined on $\Re^n$ and is continuous. 
\item[(c)]If, in addition, $f$ is nonnegative and satisfies $f(0)=0$, then $\xi(0)=0$. 
\end{enumerate}
\end{lem}
%%%%%
\begin{pf}
We first address the existence of a solution $\xi$  to the equation $\xi\in \partial f(x-W\xi)$. 
For this purpose, consider the function $f^w:\Re^n\rightarrow\Re$ defined by $f^w(x)=f(W^{1/2}x)$. Since $f$ is proper, closed and convex, then so is $f^w$ thanks to the assumed invertibility of $W$.  By \cite[Thm 6.3]{Beck17-Book}, we know then that the function $z \mapsto \prox_{f^w}(z)$ is well-defined on $\Re^n$. It follows that for any $z\in \Re^n$, the vector $z^\star$ defined by  
$z^\star= \argmin_{u\in\Re^n}F(u,z)$, with $F(u,z)=\frac{1}{2}\left\|u-z\right\|_2^2+f^w(u)$, 
exists and is unique.  Let  $x = W^{1/2}z$ and $x^\star = W^{1/2}z^\star$.  
The optimality condition of the prox operator reads as $0\in \partial_u F(z^\star,z)=z^\star-z+W^{1/2}\partial f(W^{1/2}z^\star)$, which means that there exists a vector $\xi\in \partial f(W^{1/2}z^\star)=\partial f(x^\star)$  such that $0=z^\star-z+W^{1/2}\xi$. 
Multiplying on the left by $W^{1/2}$ and re-arranging gives  $x^\star =x-W\xi^\star$ with $\xi^\star\in \partial f(x^\star)=\partial f(x-W\xi^\star)$. Since the existence and uniqueness of $x^\star$ are guaranteed, the same holds for $\xi^\star$, which finishes the proof of statement (a).\\
%%%%
Consider now proving the continuity statement (b). For this purpose, let $\left\{x_n\right\}$ be a sequence converging to $x$ as $n\rightarrow+\infty$. 
We just need to show that $\left\{\xi(x_n)\right\}$ then converges to $\xi(x)$. Using the subdifferential property as above, we can write
$$\begin{aligned}
	f(x-W\xi(x))-f(x_n-W\xi(x_n)) \geq &\xi(x_n)^\top \big[x-x_n\big. \\
	&\big.-W(\xi(x)-\xi(x_n))\big],
\end{aligned}$$
and
$$\begin{aligned}
	f(x_n-W\xi(x_n)) -f(x-W\xi(x))\geq &\xi(x)^\top \big[x_n-x\big. \\
	&\big.-W(\xi(x_n)-\xi(x))\big].
\end{aligned}$$
Summing these inequalities give
$$\begin{aligned}
	0\geq &\big(\xi(x)-\xi(x_n)\big)^\top (x_n-x)\\
	&\qquad \qquad  +\big(\xi(x)-\xi(x_n)\big)^\top W \big(\xi(x)-\xi(x_n)\big). 
\end{aligned}$$
From this, it is immediate that the sequence $\left\{\xi(x_n)\right\}$ converges to $\xi(x)$ hence proving continuity of $x\mapsto\xi(x)$. \\
\textit{Proof of (c):} Note that $\xi(0)$ satisfies $\xi(0)\in \partial f(-W\xi(0))$. 
Consequently, $f(z)-f(-W\xi(0))\geq \xi(0)^\top (z+W\xi(0))$ $\forall z\in \Re^n$. Taking $z=0$ and using the fact that $f$ is nonnegative,  $f(0)=0$ and $W$ is positive definite, gives the result. 
\qed
\end{pf}
\noindent Consider the following assumptions, where the notations $D_t(e)$ and $\Sigma_t(e)$ are defined as in \eqref{eq:Dt} and \eqref{eq:Sigmat}, respectively:

\noindent \textbf{Assumptions. }
\begin{enumerate}[label=\textbf{A.\arabic *}]
\item The loss function $\psi$ in \eqref{eq:gt} is convex, continuous and positive definite. \label{A:positive-definite} 
\item  The weighting matrices $\left\{W_t,V_t\right\}$ are symmetric, positive-definite  and  bounded in the sense that there exist positive real numbers $(a_w,b_w)$  and $(a_v,b_v)$ such that $a_w I\preceq W_t \preceq {b}_w I$ and $a_v I\preceq V_t \preceq b_v I$. \label{A:weighting-matrices}
	\item  The system matrices $\left\{A_t,C_t\right\}$ are bounded and satisfy, together with the weighting matrices,  \\
	$D_t(e)\leq 0$ $\forall (t,e)\in \mathbb{N}_{\geq 1}\times \Re^n$. \label{A:Dt-negatif}
	\item There exists a fixed $T$ and a nondecreasing and positive definite function $\alpha_0:\Re\rightarrow\Re_+$ such that $\Sigma_t(e)$ defined in \eqref{eq:Sigmat} satisfies 
	$\Sigma_t(e)\geq \alpha_0(\left\|e\right\|)$ $\forall (t,e)\in \mathbb{N}_{\geq T}\times \Re^n$. \label{A:obsv}
\end{enumerate}
The condition $D_t(e)\leq 0$ in Assumption \ref{A:Dt-negatif} is reminiscent of the  Lyapunov equation.  
As to Assumption \ref{A:obsv} we will show later in Lemma \ref{def:uniform-observability} that it is indeed related to the observability of the system.  
\begin{thm}\label{thm-Stability}
Consider the estimator \eqref{eq:prediction-update}-\eqref{eq:measurement-update} for system \eqref{eq:system} where $f_t$ and $g_t$ are defined as in \eqref{eq:linear-ft} and \eqref{eq:gt}, respectively.   
Let Assumptions  \ref{A:positive-definite}-\ref{A:obsv} hold. 
Assume that  the uncertainties $w_t$ and $v_t$ are identically equal to zero. 
Then 
the state estimation error $e_t$ generated by the estimator \eqref{eq:prediction-update}-\eqref{eq:measurement-update} converges to $0$ as $t\rightarrow+\infty$. 
\end{thm}
%%%%
\begin{pf}
By applying the change of variable $z\rightarrow x_t+e$ in \eqref{eq:measurement-update}, we can immediately observe that the estimation error is given by: 
$$e_{t}=\argmin_{e\in \Re^n}\left[\dfrac{1}{2}\big\|W_t^{-1}\left(e-	e_{t|t-1}\right)\big\|_2^2+g_t(x_t+e)\right],$$
where $e_{t|t-1}=\hat{x}_{t|t-1}-x_t = A_{t-1}e_{t-1}$.  
Writing the optimality condition gives 
\begin{equation}\label{eq:error-dynamics}
	e_{t}  =A_{t-1}e_{t-1}-W_t^{2}\xi_t(e_t), \quad \xi_t(e_t)\in \partial g_t(x_t+e_t). 
\end{equation}
For simplicity, we will use the notation $\xi_t$ instead of $\xi_t(e_t)$. 
Consider now the functions $F_t:\Re^n\rightarrow\Re_+$ defined by  $F_t(e)=e^\top W_t^{-2}e$. 
We have 
\begin{equation}\label{eq:Delta-V}
	\begin{aligned}
		F_t(e_t)-F_{t-1}(e_{t-1}) =&e_{t-1}^\top \left(A_{t-1}^\top W_t^{-2}A_{t-1}-W_{t-1}^{-2}\right)e_{t-1}\\
		&-2e_{t-1}^\top A_{t-1}^\top \xi_t+\xi_t^\top W_t^{2}\xi_t.
	\end{aligned}
\end{equation}
Since $\xi_t\in \partial g_t(x_t+e_t)=\partial \tilde{g}_t(e_t)$ (in the assumed noise-free setting), we have
\begin{equation}\label{eq:subgradient-property}
	g_t(z)-g_t(x_t+e_t)\geq \xi_t^\top \left(z-x_t-e_t\right) \quad \forall z \in \Re^n.
\end{equation}
In particular, by letting $z=x_t$, we get
$$g_t(x_t)-g_t(x_t+e_t)\geq -\xi_t^\top e_t =-\xi_t^\top A_{t-1}e_{t-1}+\xi_t^\top W_t^{2}\xi_t $$
and hence,
$$ -2\xi_t^\top A_{t-1}e_{t-1}\leq 2\left[g_t(x_t)-g_t(x_t+e_t)\right]-2\xi_t^\top W_t^{2}\xi_t. $$
Combining this with \eqref{eq:Delta-V} yields 
$$\begin{aligned}
	F_t(e_t)-F_{t-1}(e_{t-1}) &\leq e_{t-1}^\top \left(A_{t-1}^\top W_t^{-2}A_{t-1}-W_{t-1}^{-2}\right)e_{t-1}\\
	& +2\left[g_t(x_t)-g_t(x_t+e_t)\right]-\xi_t^\top W_t^{2}\xi_t. 
\end{aligned}$$
Note that $g_t(x_t)=0$ and  $g_t(x_t+e_t)=\psi(V_t^{-1}C_te_t)$ when $v_t=0$. 
Therefore by recalling the definition of $G_t$ in \eqref{eq:Gt}, we can write
\begin{equation}\label{eq:Gt-decreasing}
	G_t(e_t)-G_{t-1}(e_{t-1})\leq D_t(e_{t-1})-\xi_{t-1}^\top W_{t-1}^{2}\xi_{t-1}.
\end{equation}
Under Assumption \ref{A:Dt-negatif}, we obtain
\begin{equation}\label{eq:decreasing-G_t-I}
G_t(e_t)-G_{t-1}(e_{t-1})\leq -\xi_{t-1}^\top W_{t-1}^{2}\xi_{t-1}.
\end{equation}
This indicates that $\left\{G_t(e_t)\right\}$ is a nonnegative and nonincreasing sequence, and therefore, it is upper-bounded and convergent. Further, the above inequality implies that 
\begin{equation}\label{eq:decreasing-G_t-II}
G_t(e_t)-G_{t-T}(e_{t-T})\leq -\sum_{k=t-T}^{t-1}\xi_{k}^\top W_{k}^{2}\xi_{k}. 
\end{equation}
By Assumption \ref{A:obsv}, it follows that $G_t(e_t)-G_{t-T}(e_{t-T})\leq -\alpha_0(\left\|e_{t-T}\right\|)$. 
Note further that the inequality still holds if we replace $T$ with any $T'\geq T$. 
From this we assert that $\left\|e_t\right\|\rightarrow 0$ when $t\rightarrow 0$. 
Suppose by contradiction that this is not the case. 
Then there exists $\delta>0$  and an increasing  sequence $\{k_j\}_{j\in \mathbb{N}}$ of times in $\mathbb{N}$  such that $k_j\to \infty$ as $j\to \infty$, $k_j-k_{j-1}\geq T$ for $j\geq 1$ and $\left\|e_{k_j}\right\|>\delta$ for all $j$. Because $\alpha_0$ is nondecreasing and positive-definite, this implies that $\alpha_0(\|e_{k_j}\|)\geq \alpha_0(\delta)>0$ for all $j$.
As a result, $$G_{k_{j+1}}(e_{k_{j+1}})\leq G_{k_j}(e_{k_j}) -\alpha_0(\delta) \qquad \forall j\in \mathbb{N}.$$ 
Applying this recursively gives 
$$G_{k_j}(e_{k_j})\leq G_{k_0}(e_{k_0}) -j\alpha_0(\delta), $$
which suggests that for a sufficiently large $j$, the quantity $G_{k_j}(e_{k_j})$ will become negative, a contradiction with the fact that the function $G_t$ is nonnegative by definition. We therefore conclude that the estimation error $e_t$ must converge to $0$. \qed
\end{pf}

 %%%%%%%%%%%%%%%%%%%%%%%%%%%%%%
We will now proceed to showing that Assumption \ref{A:obsv} holds provided that the system is uniformly completely observable (UCO) in the following sense.
\begin{definition}[Generalized UCO]\label{def:uniform-observability}
Let $\psi$ and $V_k$ be subject to Assumptions  \ref{A:positive-definite} and \ref{A:weighting-matrices} respectively. 
The linear system with matrices $(A_t,C_t)$ is \textit{uniformly completely observable} if there exist class $\mathcal{K}_\infty$ functions $\underline{\alpha}$ and $\overline{\alpha}$ and a finite time horizon $T$ such that for all $e\in \Re^n$, the sequence $\left\{e_k\right\}$ defined by $e_t=e$ and $e_{k}=A_{k-1}e_{k-1}, k=t+1,\ldots,t+T$ satisfies   
\begin{equation}
	\underline{\alpha}(\left\|e\right\|) \leq \sum_{k=t}^{t+T}\psi(V_k^{-1}C_ke_k) \leq \overline{\alpha}(\left\|e\right\|) \quad \forall t\in \mathbb{N}.
\end{equation}
\end{definition}
Definition \ref{def:uniform-observability}  reduces to the standard concept of UCO \cite[Chap. 7]{Jazwinski70-Book} if $\psi$ is quadratic of the form $\psi(y)=y^\top y$ and $\underline{\alpha}(r)=\underline{q}r^2$ and $\overline{\alpha}(r)=\overline{q}r^2$ for some positive numbers $\underline{q}$ and $\overline{q}$. 
%%%%%%%%%  
\begin{lem}\label{lem:obsv-equivalence}
Let Assumptions \ref{A:positive-definite} and \ref{A:weighting-matrices} hold and consider a noise-free setting. 
If  the LTV system with matrices $(A_t,C_t)$ is uniformly completely observable in the sense of Definition \ref{def:uniform-observability},  then Assumption \ref{A:obsv} is satisfied. 
\end{lem}
%%%%
\begin{pf}
By Lemma \ref{lem:uniqueness}, part (c), it can be seen from \eqref{eq:Sigmat} that $\Sigma_t(0)=0$ in the noise-free setting. Assume that  $\Sigma_t(e)=0$. Then it follows from \eqref{eq:Sigmat} that $\xi_k(e_k)=0$ for all $k=t-T,\ldots,t-1$ because the matrices $W_k^2$ are positive definite by Assumption  \ref{A:weighting-matrices}. Hence $e_k=A_{k-1}e_{k-1}$ for $k=t-T+1,\ldots,t-1$ and $e_{t-T}=e$.  
On the other hand $\xi_k(e_k)=0\in \partial \tilde{g}_k(e_k)$ implies that for all $k=t-T,\ldots,t-1$, 
$e_k$ is a minimizer of $\tilde{g}_k$.  However the minimum of $\tilde{g}_k$ is clearly equal to zero. 
Hence $\tilde{g}_k(e_k)=\psi(V_k^{-1}C_ke_k)=0$ for all $k=t-T,\ldots,t-1$. 
By the property of uniform observability of $(A_t,C_t)$, it must hold that  $e=0$. We therefore conclude that $\Sigma_t$ is positive definite for any $t\geq T$. 
Consider now the function $\sigma:\Re_+\rightarrow\Re_+$ defined by 
$\sigma(r) = \inf_{(t,e):\left\|e\right\|\geq r, t\geq T}\Sigma_t(e)$.  
Then $\sigma$ is easily seen to be nondecreasing. Moreover, $\sigma(r)\geq 0$ and $\sigma(0)=0$. Assume for contradiction purpose that $\sigma(r)=0$ for some $r>0$. 
Then for any $\epsilon>$  there exist $(t,e)$ with $t\geq T$ and $\left\|e\right\|\geq r$ such that $\Sigma_{t}(e)<\epsilon$. This, by the boundedness properties of $W_k$ (Assumption  \ref{A:weighting-matrices}), implies that $a_w^2\sum_{k=t-T}^{t-1}\xi_k^\top \xi_k<\epsilon$ and hence that $\left\|\xi_k\right\|_2<\sqrt{\epsilon}/a_w$ for $k=t-T,\ldots,t-1$. 
On the other hand, since $\xi_k\in \partial \tilde{g}_k(e_k)$, we have 
$\tilde{g}_k(z)-\tilde{g}_k(e_k)\geq \xi_k^\top (z-e_k)$ for all $z\in \Re^n$.
In particular, for $z=0$, we get $0\leq \tilde{g}_k(e_k)=\psi(V_k^{-1}C_ke_k)\leq \xi_k^\top e_k\leq \left\|e_k\right\|_2\sqrt{\epsilon}/a_w$. 
Taking the sum gives $\sum_{k=t-T}^{t-1}\psi(V_k^{-1}C_ke_k)\leq \sqrt{\epsilon}/a_wK(e)$ with $K(e)=\sum_{k=t-T}^{t-1}\left\|e_k\right\|_2\leq K_1\left\|e\right\|+K_2b_w^2\sqrt{\epsilon}/a_w$ for some numbers $K_1$ and $K_2$ which are independent of $\epsilon$.
Hence we can make  $\sum_{k=t-T}^{t-1}\psi(V_k^{-1}C_ke_k)$  arbitrarily small for a sequence obeying \eqref{eq:Sigmat} with arbitrary small $\xi_k$ and initial error $e$ such that $\left\|e\right\|\geq r>0$. This contradicts the uniform observability condition. We therefore conclude that $\sigma$ is a positive-definite function and so, Assumption \ref{A:obsv} holds since $\Sigma_t(e)\geq \sigma(\|e\|)$ for all $t\geq T$. \qed
\end{pf}
\color{black}

\subsection{The Kalman filter as an illustrative example}\label{subsec:kalman}
For illustration purposes, we can interpret the current estimation framework in the special case of the Kalman filter (even though this paper does not  use stochastic modeling of uncertainties). 
We do so by relying on similar notations as in \cite{Zhang17-SCL}. Let the loss function $\psi$ in \eqref{eq:gt} be quadratic of the form $\psi_{\qd}(y)=(y^\top y)/2$. We select the weighting matrices $V_t$ and $W_t$ such that $V_t^2 = R_t$, where $R_t$ is a surrogate of \textit{the measurement noise covariance},  and 
$$W_t^2=P_{t|t-1}=A_{t-1}P_{t-1|t-1}A_{t-1}^\top +Q_{t-1},$$
 with $Q_{t-1}$ playing the role of  \textit{the covariance matrix of the process noise}, $w_{t-1}$ and $P_{t-1|t-1}$  defined by 
\begin{equation}\label{eq:Pt-1|t-1}
	P_{t-1|t-1}=\big(W_{t-1}^{-2}+C_{t-1}^\top V_{t-1}^{-2} C_{t-1}\big)^{-1}.
\end{equation}
Then the estimator \eqref{eq:prediction-update}-\eqref{eq:measurement-update} is equivalent to the Kalman filter for linear dynamics, as in \eqref{eq:linear-ft}. 
In this setting, assume that the matrices $Q_t$, $R_t$, $A_t$, $C_t$ are uniformly bounded in the following sense: there exist  constant and positive numbers $q_i,r_i, c_i,a_i$, $i\in \left\{1,2\right\}$, such that 
$q_1I_n\preceq Q_t\preceq q_2I_n$, $r_1I_{n_y}\preceq R_t\preceq r_2I_{n_y}$, $c_1I_{n_y}\preceq C_tC_t^\top\preceq c_2I_{n_y}$, $A_tA_t^\top \preceq a_2 I_n$. Then the matrices $P_{t|t-1}$ and $P_{t|t}$ are uniformly bounded and invertible provided that  $P_0\succ 0$, $A_t$ is nonsingular and $(A_t,C_t)$ is UCO, see Definition \ref{def:uniform-observability} above and \cite[Chap. 7]{Jazwinski70-Book}. It can then be shown that Assumption \ref{A:Dt-negatif} is satisfied in this case. 
To see this, let $D_t^{\qd}$ denote the function in \eqref{eq:Dt} when it is based on the quadratic loss $\psi_{\qd}$. Then 
$$
\begin{aligned}
&D_t^{\qd}(e) \\ 
& =\!e^\top \left[A_{t-1}^\top P_{t|t-1}^{-1}A_{t-1}-P_{t-1|t-2}^{-1}-C_{t-1}^\top V_{t-1}^{-2}C_{t-1}\right]e\\
&= \!e^\top\! \! \left[\! A_{t-1}^\top \!\big(A_{t-1}P_{t-1|t-1}A_{t-1}^\top \!\!+\!Q_{t-1}\big)^{\!-1}\!A_{t-1}\!-\!P_{t-1|t-1}^{-1}\!\right]\!e\\
& =\! -e^\top\! \Big[P_{t-1|t-1}^{-1}\!\left( P_{t-1|t-1}^{-1}\!+\!A_{t-1}^\top Q_{t-1}^{-1}A_{t-1}\right)^{\!-1}\!\! P_{t-1|t-1}^{-1}\Big]e 
\end{aligned}
$$
where we used the Woodbury matrix identity in the last line. Hence $D_t^{\qd}(e)$ is nonpositive for all $e\in \Re^n$.

Now, let us check Assumption \ref{A:obsv}. 
When $\psi$ is quadratic as given above, we get $\xi_k(e_k)=C_k^\top R_k^{-1}C_k e_k$, so that the error dynamics takes the form $M_ke_k=A_{k-1}e_{k-1}$ with $M_k=I_n+P_{k|k-1}C_k^\top R_k^{-1}C_k$. Then it follows from \eqref{eq:Sigmat} that
$$\Sigma_t(e)=\sum_{k=t-T}^{t-1}e_k^\top C_k^\top \Omega_kC_k e_k\succeq \omega_0 \sum_{k=t-T}^{t-1}e_k^\top C_k^\top C_k e_k$$
subject to $e_k=\bar{A}_{k-1}e_{k-1}$, $e_{t-T}=e$, with $\bar{A}_{k-1}=M_k^{-1}A_{k-1}$ and $\Omega_k=R_k^{-1}C_kP_{k|k-1}C_k^\top R_k^{-1}$. The rightmost inequality  follows from the assumptions above with $\omega_0=q_1c_1r_2^{-2}$. 
From this we see that  Assumption  \ref{A:obsv} is fulfilled if the pair $\left\{\bar{A}_k,C_k\right\}$ is UCO. Such an observability condition, however, concerns the closed-loop system. Interestingly, it can be shown that if the open loop system $\left\{A_k,C_k\right\}$ is UCO, then so is the closed-loop $\left\{\bar{A}_k,C_k\right\}$. This fact is proved in Section \ref{subsec:UCO} and  Lemma \ref{lem:equivalence-UCO} of the appendix. 

In conclusion, the Kalman filter can be viewed as a special case of our framework corresponding to a particular choice of the weighting matrices $W_t$ and of the loss function $\psi$. In this case, given $W_0^2\succ 0$ and the sequences $\{Q_t\}$, $\{V_t\}$ all symmetric and positive definite, the weighting matrices $\{W_t\}$ obey the Riccati equation
\begin{equation}\label{eq:Wt-Seq}
\begin{aligned}
\!\!W_t^2&=A_{t-1}W_{t-1}^2A_{t-1}^\top +Q_{t-1}\\
	&\hspace{-0.7cm}  -\!A_{t-1}W_{t-1}^2C_{t-1}^\top\!\left(V_{t-1}^2\!+\!C_{t-1}W_{t-1}^2C_{t-1}^\top\right)^{\!-1}\!C_{t-1}W_{t-1}^2A_{t-1}^\top.
	\end{aligned}
\end{equation}

\subsection{Further discussions of Assumption \ref{A:Dt-negatif}}
In this section we illustrate two possible ways of selecting the weighting matrices $\left\{W_t\right\}$, $\left\{V_t\right\}$ entering the sequence of optimization problems \eqref{eq:measurement-update} so that Assumption \ref{A:Dt-negatif} holds.

A first observation is that if\footnote{This condition accounts for the situation where $\sum_{t=0}^{\infty}\log \left\|A_t\right\|_2=-\infty$} $\sum_{t=0}^{\infty}\log \left\|A_t\right\|_2<+\infty$, then for any given sequence $\left\{V_t\right\}$ obeying Assumption \ref{A:weighting-matrices}, we can easily select a bounded sequence  $\left\{W_t\right\}$ such that Assumption \ref{A:Dt-negatif} holds true. It suffices for this that $A_{t-1}^\top W_{t}^{-2}A_{t-1}-W_{t-1}^{-2}\preceq 0$ or, equivalently, that $\big\|W_{t}^{-1}A_{t-1}W_{t-1}\big\|_2\leq 1$, which yields $\gamma_t\geq \gamma_{t-1}\left\|A_{t-1}\right\|_2$ if we choose $W_k$ in the form $W_k=\gamma_k I$.

A second observation is that we can indeed compute $W_t$ recursively using the same update equation as in the Kalman filter (see, Eq. \eqref{eq:Wt-Seq}) even though a different loss function $\psi$ than the quadratic one is used in \eqref{eq:measurement-update}-\eqref{eq:gt}. 

\begin{lem}
Assume that the weighting matrices $W_t$ are updated according to \eqref{eq:Wt-Seq} with $\left\{V_t\right\}$ obeying Assumption \ref{A:weighting-matrices} and $\left\{Q_{t}\right\}$ such that 
\begin{equation}\label{eq:cond-for-A3}
	A_{t-1}^\top S_{t-1}^{-1}A_{t-1}-W_{t-1}^{\!-2}\preceq 0 \quad \forall t\geq 1,  
\end{equation}
with $S_{t-1}=Q_{t-1}+A_{t-1}P_{t-1|t-1}A_{t-1}^\top$ and $P_{t-1|t-1}$ as in \eqref{eq:Pt-1|t-1}. 
Then the function $D_t$ defined in \eqref{eq:Dt} satisfies Assumption  \ref{A:Dt-negatif}. 
In particular, if $A_{t-1}$ is invertible and $Q_{t-1}$ is chosen in the form  $Q_{t-1}=q_{t-1}A_{t-1}A_{t-1}^\top$, $q_{t-1}>0$, then a sufficient condition for Assumption \ref{A:Dt-negatif} to hold is that $q_{t-1}\geq \left\| W_{t-1}^2C_{t-1}^\top\!\left(V_{t-1}^2\!+\!C_{t-1}W_{t-1}^2C_{t-1}^\top\right)^{\!-1}\!C_{t-1}W_{t-1}^2\right\|_2$. 
\end{lem}
%%%
\begin{pf}
Observe that $D_t(e)$ can be expressed in terms of $D_t^{\qd}(e)$ as follows
$$
\begin{aligned}
	D_t(e) & =D_t^{\qd}(e)+2\psi_{\qd}(V_{t-1}^{-1}C_{t-1}e)-2\psi\big(V_{t-1}^{-1}C_{t-1}e\big) \\
			& = e^\top\! \left(A_{t-1}^\top S_{t-1}^{-1}A_{t-1}-W_{t-1}^{\!-2}\right)e-2\psi\big(V_{t-1}^{-1}C_{t-1}e\big)
\end{aligned}
$$
It follows that a sufficient condition for $D_t(e)\leq 0$, as required by Assumption  \ref{A:Dt-negatif}, is that 
$A_{t-1}^\top S_{t-1}^{-1}A_{t-1}-W_{t-1}^{-2}\preceq 0$ as stated in the lemma. 
With $Q_{t-1}=q_{t-1}A_{t-1}A_{t-1}^\top$, this condition becomes
$$\begin{aligned}
	&A_{t-1}^\top \Big(A_{t-1}\big(q_{t-1}I+P_{t-1|t-1}\big)A_{t-1}^\top\Big)^{\!-1}\!A_{t-1}-W_{t-1}^{-2}\preceq 0 \\
	& \Leftrightarrow \Big(q_{t-1}I+P_{t-1|t-1}\Big)^{\!-1}-W_{t-1}^{-2}\preceq 0\\
	& \Leftrightarrow  q_{t-1}I\succeq W_{t-1}^{2}-P_{t-1|t-1}.
\end{aligned}$$
Now using the Woodbury identity to develop $P_{t-1|t-1}$ from \eqref{eq:Pt-1|t-1} completes the proof of the second statement. 
\end{pf}
Note that \eqref{eq:cond-for-A3} is always achievable provided that $Q_{t-1}$ is chosen sufficiently large. 
Also, even if $A_{t-1}$ is not invertible, it is still possible, though at the price of a slightly more involved lower bound,  to pick $q_{t-1}$ such that $Q_{t-1}=q_{t-1}A_{t-1}A_{t-1}^\top$ satisfies \eqref{eq:cond-for-A3}.

\subsection{Component-wise measurement update}
For computational reasons, it may be useful to approximate, at each time $t$, the measurement update scheme by a finite sequence of scalar-valued measurement updates.  As will be substantiated in Section \ref{sec:special-cases}, this allows for obtaining closed-form expressions of the estimator.  
To present the principle of this component-wise update, assume that whenever  $V_t$ in \eqref{eq:gt} is diagonal with  $V_t^{-1}=\diag(\lambda_{t1},\cdots,\lambda_{tn_y})$ with strictly positive reals $\lambda_{ti}$, the function $g_t$ is decomposable as 
$g_t(z)= \sum_{i=1}^{n_y}g_{ti}(z) $
 with $g_{ti}(z)=\psi_0(\lambda_{ti}(y_{ti}-c_{ti}^\top z))$,  $c_{ti}^\top$ standing for the $i$-th row of $C_t$, $y_{ti}$ for the $i$-th entry of $y_t$, and $\psi_0:\Re\rightarrow\Re_+$. 
We can then consider extracting the information from the measurement vector $y_t\in \Re^{n_y}$ by processing the components of $y_t$ one after another. This gives the following iterative update equation: Set $z_t^{(0)} = f_{t-1}(\hat{x}_{t-1})$ and for $i=1,\ldots, n_y$, let  $z_t^{(i)}$ be defined by 
\begin{equation}\label{eq:componentwise}
%\left\{	\begin{aligned}
	%	& z_t^{(0)} = f_{t-1}(\hat{x}_{t-1}) \\ 
		z_t^{(i)} = \argmin_{z\in \Re^n}\left[\dfrac{1}{2}\big\|W_{t}^{-1}(z-z_t^{(i-1)})\big\|_2^2+g_{ti}(z)\right]
	%	& \qquad \qquad \qquad \qquad i=1,\ldots, n_y\\
	% & \hat{x}_{t} =  z_t^{(n_y)}
	%\end{aligned}\right.
\end{equation}
and then  set the estimate at time $t$ to be  $\hat{x}_{t} =  z_t^{(n_y)}$, the $n_y$-th iterate. 
Convergence of this (relaxed) estimator can be studied under similar conditions as in Theorem \ref{thm-Stability}. To proceed, let us  change the time scale by introducing a new time variable $\tau$ which  concatenates the true discrete time $t$ and  the $n_y$ measurement update iterations at each $t$. That is, consider the mapping $\tau:\mathbb{N}\times \left\{0,1,\ldots,n_y\right\}\rightarrow \mathbb{N}$, $(t,i)\mapsto \tau=\tau(t,i)=tn_y+i$ with the constraint that $\tau(t,n_y)=\tau(t+1,0)$. 

Then, by replacing $e_t^{(i)}\triangleq z_t^{(i)}-x_t$ with $\bar{e}_\tau$, the  estimation error dynamics \eqref{eq:error-dynamics} can be written in the form of a \textit{single output} time-varying system 
\begin{equation}\label{eq:ebar}
	\bar{e}_{\tau} = \hat{A}_{\tau-1}\bar{e}_{\tau-1} -W_\tau^2\bar{\xi}_{\tau}(\bar{e}_{\tau}), 
\end{equation}
where $\hat{A}_{\tau-1}=A_{t-1}$ if $\tau=tn_y+1$ and $\hat{A}_{\tau-1}=I$ otherwise; $\bar{\xi}_{\tau}\in \partial\psi_0(\lambda_{\tau}\bar{c}_{\tau}^\top \bar{e}_{\tau})$ with $\bar{c}_{\tau}^\top=c_{ti}^\top$ for $\tau=tn_y+i$, $i=1,\ldots,n_y$. 

\begin{thm}[Convergence of the component-wise update]
Consider the estimator \eqref{eq:componentwise} with linear system dynamics as described in \eqref{eq:linear-ft}. Let Assumptions \ref{A:positive-definite}-\ref{A:obsv} hold with $\psi$ replaced by $\psi_0$ and $D_t(e)$ replaced by $\bar{D}_t(\bar{e})$ defined by  
$$
\begin{aligned}
	\!\bar{D}_t(\bar{e})\!= \!\!\sum_{\tau=tn_y+1}^{(t+1)n_y}&\!\left[\bar{e}_{\tau-1}^\top \left(\hat{A}_{\tau-1}^\top W_{\tau}^{-2}\hat{A}_{\tau-1}-W_{\tau-1}^{-2}\right)\bar{e}_{\tau-1}\right. \\
										& \hspace{1cm} \left.-2\psi_0(\lambda_{\tau-1}\bar{c}_{\tau-1}^\top \bar{e}_{\tau-1})\right] 
	\end{aligned}
$$
subject to \eqref{eq:ebar} and $\bar{e}=\bar{e}_{tn_y}$.  
If the  functions  $g_{ti}$ are proper, closed, convex, and the uncertainties $w_t$ and $v_t$ are identically equal to zero, 
then 
the state estimation error $\left\{\bar{e}_{tn_y}\right\}$ generated by the estimator \eqref{eq:componentwise} converges to $0$ as $t\rightarrow+\infty$. 
\end{thm}
%%%%%%%%%%
%\color{blue}
\subsection{Robustness to impulsive measurement noise}
In the case where the open-loop system is uniformly exponentially stable, we can show that the estimation error  \eqref{eq:error-dynamics-linear} induced by the proximal observer is bounded 
provided that the loss function $\psi$ has bounded subgradients. An important feature of this result is that the resulting  bound is independent of any bound on the measurement noise $\left\{v_t\right\}$. It follows that the estimator is robust to arbitrarily large measurement noise. 
%%%
\begin{prop}
Consider the estimator \eqref{eq:prediction-update}-\eqref{eq:measurement-update} for system \eqref{eq:system} where $f_t$ and $g_t$ are as in \eqref{eq:linear-ft} and \eqref{eq:gt} respectively; the process noise $\left\{w_t\right\}$ is bounded while the measurement noise can take on possibly arbitrarily large values.   
Let Assumption \ref{A:weighting-matrices} hold and  $(A_t,C_t)$ be bounded.   
Suppose further that the system defined by $z_t=A_{t-1}z_{t-1}$ is uniformly exponentially stable in the sense that there exists $c>0$ and $\lambda\in (0,1)$ such that its transition matrix $\Phi$ obeys $\left\|\Phi(t,t_0)\right\|\leq c\lambda^{t-t_0}$ $\forall t\geq t_0$. 
Under these conditions, if $\psi$ satisfies the triangle inequality property, then the  (estimation) error system \eqref{eq:error-dynamics-linear} is bounded as follows 
\begin{multline}
\|e_t\|_2\leq \beta(\|e_0\|_2,t) + \gamma_w(\max_{s\leq t-1}\{\|w_s\|_2\}) + R \ , \forall t\geq0
\end{multline}
where the functions $\beta$,  $\gamma_w$ are defined by 
$\beta(r,t)=cr\lambda^t$, $\gamma_w(r)=qr$ and  the constant $R$ is given by $R=q\delta_\psi\sqrt{n}\sup_{t}\|W_t^2C_t^\top V_t^{-1}\|_2$, with $q=c/(1-\lambda)$ and $\delta_\psi$ being a constant depending solely on the loss function $\psi$. 
\end{prop}
\begin{pf}
Noting that $\xi_t(e_t)=-C_t^\top V_t^{-1}h_t(e_t)$, with $h_t(e_t)\in \partial\psi\big(V_t^{-1}(v_t-C_te_t)\big)$, we can rewrite the error dynamics in \eqref{eq:error-dynamics-linear} as 
$e_t = A_{t-1}e_{t-1}+\bar{w}_{t-1}$ 
where $\bar{w}_{t-1}=W_t^2C_t^\top V_t^{-1}h_t(e_t)-w_{t-1}$. \\
\textit{Claim:} There exists  a constant bound $\delta_{\psi}$ such that  $\left\|h_t(e)\right\|_\infty\leq \delta_{\psi}$ for all $(t,e)$. To see this, let us show that   $\left\|h\right\|_\infty\leq \delta_{\psi}$ for all $(h,y)\in \Re^{n_y}\times \Re^{n_y}$ related by $h\in \partial \psi(y)$. 
Indeed, for such a pair $(h,y)$, we have $h^\top y'\leq \psi(y+y')-\psi(y)$ for all $y'\in \Re^{n_y}$. This implies, by the triangle inequality property of $\psi$,  that $h^\top y'\leq \psi(y')$. In particular, taking $y'$ equal to the canonical vectors $\pm p_i$ of $\Re^{n_y}$, we get $|h_i|\leq \max(\psi(p_i),\psi(-p_i))$, $i=1,\ldots,n_y$. This shows that $\left\|h\right\|_\infty\leq \delta_{\psi}$ with the bound $\delta_{\psi}=\max_{i=1,\ldots,n_y}\max(\psi(p_i),\psi(-p_i))$ depending only on $\psi$. Applying this with $y=V_t^{-1}(v_t-C_te_t)$ and $h=h_t(e)$ proves the claim.  \\
Using now the other assumptions of the proposition we have that the term $W_t^2C_t^\top V_t^{-1}h_t(e_t)$ is uniformly bounded by $R/q$ hence, implying finally that $\left\{\bar{w}_{t-1}\right\}$ is bounded. The rest of the proof follows by bounding directly the norm of the  solution of the error system, $e_t=\Phi(t,0)e_0+\sum_{j=0}^{t-1}\Phi(t,j)\bar{w}_{j}$. \qed
\end{pf}
Note that the assumption of the triangle inequality property is realized for example if $\psi$ is a norm. Indeed, this assumption can be weakened by requiring that the  property only holds outside some compact set on which $\psi$ is continuous. 
\color{black}

\subsection{Convergence in the nonlinear case}\label{subsec:analysis-NL}
Consider now the more general situation where $f_t$ is nonlinear. For the purpose of analyzing the estimation error \eqref{eq:error}, we require the following assumption.
\begin{enumerate}[resume,label=\textbf{A.\arabic *}]
\item Each map $f_t$ from \eqref{eq:system} is globally Lipschitz, i.e., for all $t\in\mathbb{N}$, there exists $L_t\geq 0$ such that $\left\|f_t(x)-f_t(x')\right\|\leq L_t\left\|x-x'\right\|$ $\forall x,x'\in \Re^n$.\label{assum:Lipschitz}
\end{enumerate}
Moreover, we change the definition of $\Sigma_t$ to match the nonlinear dynamics of the estimation error as described in \eqref{eq:error}. To this end, we introduce the following notation: 
$$
\begin{aligned}
	& \widetilde{\Sigma}_t(e) =\sum_{k=t-T}^{t-1}\xi_k(e_k)^\top W_k^2\xi_k(e_k) \notag \\
							& \qquad \qquad  \mbox{subject to }  e_k=d_{k-1}(e_{k-1})-W_k^2\xi_k(e_k), \label{eq:Sigmat2}\\
							& \qquad \qquad \qquad \qquad \: \:   e_{t-T}=e \notag 
	\end{aligned}
$$
where $d_{k-1}(e_{k-1}) = f_{k-1}({x}_{k-1}+e_{k-1})-f_{k-1}(x_{k-1})$. 
In addition, let 
$$
 \widetilde{D}_t(e)=e^\top\left(a_w^{-2}L_{t-1}^2I_n-W_{t-1}^{-2}\right)e-2\psi(V_{t-1}^{-1}C_{t-1}e)
$$ 
with $L_{t}$ being the Lipschitz constant of the maps $f_t$ and $a_w$ the lower bound on $\left\{W_t\right\}$  from Assumption \ref{A:weighting-matrices}.

\begin{thm}[Convergence -- nonlinear case] $\: $\\
Consider the system \eqref{eq:system} and the estimator \eqref{eq:prediction-update}-\eqref{eq:measurement-update} under the assumptions that  
(i) the uncertainties $w_t$ and $v_t$ are identically equal to zero;  (ii) Assumptions  \ref{A:positive-definite}, \ref{A:weighting-matrices} and  \ref{assum:Lipschitz} hold true; (iii) Assumption \ref{A:obsv} is true with $\Sigma_t$ replaced by $\widetilde{\Sigma}_t$; (iv)  $\widetilde{D}_t(e)\leq 0$ $\forall (t,e)\in \mathbb{N}_{\geq 1}\times \Re^n$.   \\
Then the estimation error in \eqref{eq:error} satisfies $e_t\rightarrow 0$ as $t\rightarrow\infty$. 
\end{thm}
%%%%%
\begin{pf}
The proof of this theorem follows almost the same lines as that of Theorem \ref{thm-Stability}. 
Therefore we only give here  the key ideas. By the Lipschitz condition on $f_k$ we can write \eqref{eq:error} as $e_k=d_{k-1}-W_k^2\xi_k(e_k)$ (in the absence of noise), where $\left\|d_{k-1}\right\|\leq L_{k-1} \left\|e_{k-1}\right\|$. Note that we have used $d_{k-1}$ instead of $d_{k-1}(e_{k-1})$ for simplicity. 
Then \eqref{eq:Delta-V} becomes 
\begin{equation}\label{eq:Delta-V2}
	\begin{aligned}
		F_t(e_t)-F_{t-1}(e_{t-1}) =& \: d_{t-1}^\top W_t^{-2}d_{t-1}-e_{t-1}^\top W_{t-1}^{-2}e_{t-1}\\
		&\qquad -2d_{t-1}^\top\xi_t+\xi_t^\top W_t^{2}\xi_t\\
		&\leq e_{t-1}^\top\left(a_w^{-2}L_{t-1}^2I_n-W_{t-1}^{-2}\right)e_{t-1}\\
		& \qquad -2d_{t-1}^\top\xi_t+\xi_t^\top W_t^{2}\xi_t,
	\end{aligned}
\end{equation}
where the inequality is a consequence of the boundedness condition on $W_t$ from \ref{A:weighting-matrices}. 
The rest of the proof follows from that of Theorem \ref{thm-Stability},  essentially by replacing $A_{t-1}e_{t-1}$,  $\Sigma_t$ and $D_t$ by $d_{t-1}$,  $\widetilde{\Sigma}_t$,  and $\widetilde{D}_t$, respectively. \qed
\end{pf}

%%%%%%%%%%%%%%%%%%%%%%%%%%%%%%%%%%%%%%%%%%%%%%%%%%%%%%%%%%%%%%%%%%
A fundamental question now is how to select the loss function $\psi$ (and hence the function $g_t$) for the resulting state estimator \eqref{eq:prediction-update}-\eqref{eq:measurement-update} to enjoy robustness properties against impulsive noise.  We will study two specific choices for the function $g_t$ and show both theoretically and numerically how much the associated estimators can mitigate the adversarial effect of the noise $\left\{\zeta_t\right\}$. Three other robustness-inducing candidate loss functions are considered in the appendices. We show for all of them that by adopting a component-wise measurement update strategy (i.e., one scalar-valued measurement at a time), numerically simple closed-form update formulas can be derived. 

\section{Some special cases}\label{sec:special-cases}
This section discusses a few special cases of the framework presented in the previous section.

\subsection{An $\ell_1$ loss-based proximal observer}\label{subsec:L1}
We first consider the scenario where the loss function $\psi$ in \eqref{eq:gt} is taken to be of the form 
$\psi(e)=\left\|e\right\|_1$. In this case $g_t$ is given by
$g_t(z)=\|V_t^{-1}\left(y_t-C_t z\right)\|_{1}. $ 
Even though convex, $g_t$ is a nonsmooth (non-differentiable) function for which the state $\hat{x}_t$ in \eqref{eq:measurement-update} does not admit a closed form expression.  As discussed earlier in the paper, one can consider solving online at each time $t$ the underlying optimization problem to obtain an estimate of the state. This may however be computationally expensive for some applications as it typically  requires iterative numerical schemes.  Therefore, an appealing alternative to direct (and costly) online optimization is to consider an approximate but cheaper solution. 
We propose here a component-wise relaxation of the problem consisting in updating the prior estimate  incrementally by adding one scalar-valued  measurement (from one sensor) at a time.  To be more specific, let $y_{ti}$, $i=1,\ldots,n_y$, denote the $i$-th component of the vector $y_t$. The idea is to update, at time $t$, the prior estimate $\hat{x}_{t|t-1}$ in $n_y$ steps  using a single measurement $y_{ti}$ at each step.

\paragraph{A component-wise relaxation} 
Assume that $V_t^{-1}$ is diagonal with diagonal elements denoted  by $\lambda_{ti}$, i.e.,  $V_t^{-1}=\diag(\lambda_{t1},\cdots,\lambda_{tn_y})$. Then by applying the sensor-wise sequential update strategy, we get a new state estimator  given by 
\begin{equation}\label{eq:relaxation}
\left\{	\begin{aligned}
		&	\hat{x}_{t|t-1} = f_{t-1}(\hat{x}_{t-1})\\
		& z_t^{(0)} =  \hat{x}_{t|t-1}\\ 
		&	z_t^{(i)} \!=\! \argmin_{z\in \Re^n}\left[\!\dfrac{1}{2}\big\|W_{t}^{-1}(z-z_t^{(i-1)})\big\|_2^2+\!\lambda_{it}\left|y_{ti}-c_{ti}^\top z\right|\right]\\ 
		& \qquad \qquad \qquad \qquad i=1,\ldots, n_y\\
		& \hat{x}_{t} =  z_t^{(n_y)}
	\end{aligned}\right.
\end{equation}
where $c_{ti}^\top$ is the $i$-th row of the observation matrix $C_t$ and $|\cdot|$ refers to the absolute value. 
An important advantage of the relaxed estimator \eqref{eq:relaxation} is that it admits a much cheaper implementation cost. Indeed,  we can show that the update equation for $z_t^{(i)}$ has a closed-form expression.  This relies on the following lemma. 
%%%%%%%%
\begin{lem}\label{lem:absolute-value}
Let $g:\Re^n\rightarrow\Re$ be the function defined by 
$g(z) = \alpha|a^\top z+b| $ 
with $a\in \Re^n\setminus \left\{0\right\}$, $b\in \Re$, $\alpha>0$. 
Then 
\begin{equation}
	\prox_{g}(x) = x-\alpha\Sat_{1}\left(\dfrac{a^\top x+b}{\alpha \left\|a\right\|_2^2}\right) a 
\end{equation}
where $\Sat_1$ refers to the saturation function defined in \eqref{eq:Saturation}. 
\end{lem}
%%%
\begin{pf}
The proof is based  on the principle that 
$$z^*=\prox_{ g}(x)=\argmin_{z}\left[\frac{1}{2}\left\|z-x\right\|_2^2+ g(z)\right]$$
 if and only if  $0$ lies in the subdifferential of the function  $z\mapsto \frac{1}{2}\left\|z-x\right\|_2^2+ g(z)$ at $z^*$.  By noting that $g(z)=\alpha \psi(a^\top z+b)$, with $\psi:\Re\rightarrow\Re_+$  defined by $\psi(e)=|e|$, and using the subdifferential composition rule recalled in Section \ref{subsec:maths}, the optimality condition translates into  
$$0\in z^*-x+\partial g(z^*)=z^*-x+\alpha a \partial \psi (a^\top z^*+b).$$ 
That is,  
$z^*=x-\alpha r^* a $ with $r^*\in \partial \psi (a^\top z^*+b)$. 
Recalling that the subdifferential of $\psi$ is given by 
\begin{equation}\label{eq:subdifferential-absolute-value}
	\partial \psi(e)=\left\{
	\begin{array}{ll}
	\left\{\sign(e)\right\} & \mbox{if }e\neq 0\\
	\interval{-1}{1} & \mbox{if } e=0,
	\end{array}
	\right.
\end{equation}
 we obtain $r^*\in \partial \psi (a^\top z^*+b)$, with 
$$
 \partial \psi (a^\top z^*+b)=
	 \left\{
	\begin{array}{ll}
	\left\{\sign(r(x)- r^*)\right\} & \mbox{if }r(x)- r^* \neq 0 \\
	\interval{-1}{1} & \mbox{if } r(x)- r^* = 0 
	\end{array}
	\right.
$$
and $$r(x)=\dfrac{a^\top x+b}{\alpha \left\|a\right\|_2^2}. $$
Solving this equation for $r^*$ under the specified constraints yields a unique value $r^*=\Sat_1(r(x))$, 
where   $\Sat_{\alpha}:\Re\rightarrow\Re_+$ is the saturation function defined by  
\begin{equation}\label{eq:Saturation}
	\Sat_{\alpha}(e)=\left\{
	\begin{array}{ll}
		\alpha \sign(e) & \mbox{if } \left|e\right|>\alpha\\
		e & \mbox{if } \left|e\right|\leq \alpha.
	\end{array}
	\right.
\end{equation}
Finally, inserting this in the expression of $z^*$ gives the claimed result. \qed 
\end{pf}
%%%%%%%%%%%%
By applying Lemma \ref{lem:absolute-value} together with \eqref{eq:weigthed-proxi}, we obtain an analytic expression for $z_t^{(i)}$ in \eqref{eq:relaxation} as follows 
\begin{equation}\label{eq:closed-form-Abs}
z_t^{(i)}\!=\!	z_t^{(i-1)}+\lambda_{ti}\Sat_{1}\left(\dfrac{y_{ti}-c_{ti}^\top z_t^{(i-1)}}{\lambda_{ti} \left\|W_t c_{ti}\right\|_2^2}\right) W_t^2c_{ti},
\end{equation} 
Two observations can be made in view of the expression \eqref{eq:closed-form-Abs}. First we note the appearance of the saturation function as a consequence of minimizing the absolute value of the scalar-valued error. This illustrates the good robustness-inducing property of the absolute value loss in that it implicitly prevents large measurements deviations from affecting the estimation. 
Second, at each stage $i$ of the update, the estimate is modified in the direction of $W_t^2c_{ti}$ with $c_{ti}$ being the corresponding row of $C_t$. The complexity of the measurement update step is about $\mathcal{O}(n^2n_y)$ and can be reduced to $\mathcal{O}(nn_y)$ if $W_t$ and $C_t$ are completely known off-line because then, the products $W_t^2 c_{ti}$ can be precomputed off-line and stored.

\begin{rem}
Considering the possibility of letting $\lambda_{ti}$ decrease to zero as $t\rightarrow\infty$, it is interesting to note that 
$\Sat_{\alpha}(e)=\alpha\Sat_1(e/\alpha)$ for $\alpha>0$. Hence \eqref{eq:closed-form-Abs} may alternatively be   implemented as 
$$z_t^{(i)}=	z_t^{(i-1)}+\dfrac{1}{\left\|W_t c_{ti}\right\|_2^2}\Sat_{\sigma_{ti}}\left(y_{ti}-c_{ti}^\top z_t^{(i-1)}\right) W_t^2c_{ti}, $$
with $\sigma_{ti} = \lambda_{ti}\left\|W_t c_{ti}\right\|_2^2$. 
\end{rem}

\paragraph{Stability of the relaxed estimator} 
In the case where  the system \eqref{eq:system} is an LTI single output system described by 
\begin{equation}\label{eq:LTI-Single-Output}
	\begin{aligned}
		&x_{t}=Ax_{t-1}+Bu_{t-1}, \\
		&y_t=c^\top x_t,
	\end{aligned}
\end{equation}
%If we let the weighting matrix sequence $W_t$ be constant ($W_t=W$ for all $t$), then the observer \eqref{eq:relaxation} has the form
the  observer \eqref{eq:relaxation} takes  the form
\begin{equation}\label{eq:observer-LTI}
\left\{	\begin{aligned}
		&\hat{x}_{t|t-1}=	A\hat{x}_{t-1}+Bu_{t-1}\\
		&\hat{x}_t=	\hat{x}_{t|t-1}+
		 \lambda\Sat_{1}\left(\dfrac{y_{t}-c^\top \hat{x}_{t|t-1}}{\lambda\left\|W_t c\right\|_2^2}\right) W_t^2c.  
	\end{aligned}\right. 
\end{equation}
%%%
\begin{cor}[Stability in the LTI case]\label{cor:absolute-value}
Consider the single output LTI system  \eqref{eq:LTI-Single-Output} (in a noise-free estimation setting). 
If the pair $(A,c^\top)$ is observable and satisfies 
\begin{equation}\label{eq:detectability:LTI}
e^\top \left(A^\top W_t^{-2}A-W_{t-1}^{-2}\right)e-2\lambda \left|c^\top e\right|\leq 0 \: \:  \forall e\in \Re^n,
\end{equation}
where $\left\{W_t\right\}$ is subject to the boundedness assumption \ref{A:weighting-matrices}, then  
 the error $e_t=\hat{x}_t-x_t$ associated with the estimator \eqref{eq:observer-LTI} with respect to the system  \eqref{eq:LTI-Single-Output} satisfies $\lim_{t\rightarrow\infty}e_t=0$. 
\end{cor}
%%%
\begin{pf}
The idea of the proof is to check that the conditions of Theorem \ref{thm-Stability} are satisfied and then apply that theorem to conclude. We will use the implicit expressions of the estimation error which, in this case, takes the form
$$e_{t}=\argmin_{e\in \Re^n}\left[\dfrac{1}{2}\big\|W_t^{-1}\left(e-	Ae_{t-1}\right)\big\|_2^2+\lambda \left|c^\top e\right|\right]$$ 
In this case, the loss function  $\psi$ from \eqref{eq:gt} is the absolute value, i.e., $\psi(y)=|y|$ and $V_t^{-1}=\lambda>0$. Hence the conditions \ref{A:positive-definite}-\ref{A:Dt-negatif} are satisfied. As to the condition \ref{A:obsv}, it follows from Lemma \ref{lem:obsv-equivalence}. \qed
\end{pf}
%%%%
%%%%%%%%%%%%%%%%%%%%%%%%%%%%%%%%%%%%%%%%%%%%%%%%%%%%%%

\subsection{Proximal observer based on a Lasso-type of loss}
A second proximal estimator  studied in this paper is the one for which the measurement update equation \eqref{eq:measurement-update} is defined by 
\begin{equation}
	\begin{aligned}
		& (\hat{x}_{t},\hat{\varphi}_t)=\hspace{-.3cm}\argmin_{(z,\varphi)\in \Re^n\times \Re^{n_y}}\left[\dfrac{1}{2}\big\|W_t^{-1}\left(z-	\hat{x}_{t|t-1}\right)\big\|_{2}^2+g_t^{\Lasso}(z,\varphi)\right] 
	\end{aligned}
\end{equation}
with $g_t^{\Lasso}:\Re^n\times \Re^{n_y}\rightarrow \Re_+$ given by 
\begin{equation}\label{eq:gt-Lasso}
	g_t^{\Lasso}(z,\varphi)= \dfrac{1}{2}\big\|y_t-C_t z-\varphi\big\|_{V_t^{-1}}^2+\gamma_t \big\|\varphi\big\|_1,
\end{equation}
and $\gamma_t>0$ being a regularization parameter. The additional optimization variable $\varphi$ is expected, out of the optimization,  to be an estimate of the sparse component $\zeta_t$ of the noise $v_t$ affecting the output vector $y_t$. In comparison with  the estimator of Section \ref{subsec:L1},  a quadratic loss is put on  $y_t-C_t x-\varphi$ (which is intended as an estimate of the dense component $\nu_t$ of $v_t$ in \eqref{eq:decomp-vt})  while an $\ell_1$ norm penalizes $\varphi$ (which is expected to be an estimate of the sparse component $\zeta_t$ of $v_t$). Hence, the objective reflects the assumptions we put on the uncertainties of the system model \eqref{eq:system} since as is well-known,  minimizing a quadratic loss encourages dense residual  while minimizing the $\ell_1$ norm has a sparsification property.

\paragraph{Component-wise relaxation}
Now the component-wise relaxed version of the estimate is formally given by 
\begin{equation}\label{eq:relaxation-Lasso}
	\left\{\begin{aligned}
		&	z_t^{(0)} = f_{t-1}(\hat{x}_{t-1})\\
		&	\!(z_t^{(i)},\varphi_{ti}^{(i)})\!\! = \!\!\!\!\argmin_{(z,\varphi_i)\in \Re^n\times \Re}\!\!\left[\dfrac{1}{2}\big\|W_{t}^{-1}(z-z_t^{(i-1)})\big\|_2^2\!+\!g_{ti}^{\Lasso}(z,\varphi_i)\right]\\ 
		& \qquad \qquad \qquad \qquad i=1,\ldots, n_y\\
		& \hat{x}_{t} =  z_t^{(n_y)}
	\end{aligned}\right. 
\end{equation}
with $g_{ti}^{\Lasso}:\Re^n\times \Re\rightarrow \Re_+$ denoting a component-wise version of $g_{t}^{\Lasso}$ in \eqref{eq:gt-Lasso}, 
$g_{ti}^{\Lasso}(z,\varphi_i) =\frac{\lambda_{ti}}{2}\left(y_{ti}-c_{ti}^\top z-\varphi_i\right)^2+\gamma_{ti} |\varphi_i|.$
%%%%
Here again, closed-form expressions can be derived for $(z_t^{(i)},\varphi_{ti}^{(i)})$.  
For this purpose, we will rely on the following result. 
%%%%%%%%%%%%%%%%%%%%%%%%%%%%%%%%%%%%%
\begin{lem}\label{lem:lasso-type}
Let $g:\Re^n\times \Re\rightarrow\Re$ be the function defined by 
\begin{equation}\label{eq:g(z,xi)}
	g(z,\varphi)=\dfrac{\lambda}{2}\big(a^\top z+b-\varphi\big)^2+\gamma \left|\varphi\right| 
\end{equation}
with $a\in \Re^n\setminus \left\{0\right\}$, $b\in \Re$, $\lambda>0$ and $\gamma>0$. \\
Then the pair of functions $\big(x\mapsto z_g^*(x),x\mapsto\varphi_g^*(x)\big)$ defined by 
\begin{equation}\label{eq:(z,xi)}
	\big(z_g^*(x),\varphi_g^*(x)\big)= \argmin_{(z,\varphi)\in \Re^n\times \Re}\left[\dfrac{1}{2}\left\|z-x\right\|_2^2+g(z,\varphi)\right]
\end{equation}
is well-defined on $\Re^n$ and can be expressed by 
\begin{equation}\label{eq:closed-form-zg}
\left\{	\begin{aligned}
				&z_g^*(x) = x-\gamma \Sat_{1} (\rho(x))a \\ 
				&\varphi_g^*(x)  = \eta \dz_1(\rho(x)) 
	\end{aligned}\right.
\end{equation}
where $\dz_1$ is the deadzone function defined by $\dz_1(x)=x-\Sat_1(x)$ and 
\begin{equation}\label{eq:rho(x)}
	\begin{aligned}
		&	\rho(x) = \dfrac{a^\top x +b}{\eta}, \\
		&	\eta = \gamma(\frac{1}{\lambda}+\left\|a\right\|_2^2). 
	\end{aligned}
\end{equation}
\end{lem}
%%%%%
\begin{pf}
Well-definedness of $(z_g^*,\varphi_g^*)$ consists in the existence and uniqueness of $z_g^*(x)$ and $\varphi_g^*(x)$ as defined in \eqref{eq:(z,xi)}. 
Existence is guaranteed by the fact that the objective function $(z,\varphi)\mapsto \dfrac{1}{2}\left\|z-x\right\|_2^2+g(z,\varphi)$ is coercive (i.e. radially unbounded) and continuous \cite[Thm 1.9]{Rockafellar05-Book}. As to the uniqueness we will show it below by expressing uniquely the solution. 
Noting that the objective function is convex (and hence subdifferentiable), the optimality conditions of the optimization problem in \eqref{eq:(z,xi)} can be written as
$$ 
\begin{aligned}
	& z_g^*(x)-x+\lambda\left(a^\top z_g^*(x)+b-\varphi_g^*(x)\right) a =0 \\
	& -\lambda\left(a^\top z_g^*(x)+b-\varphi_g^*(x)\right)+\gamma s=0
\end{aligned}
$$
for some $s\in \partial_{\varphi}\psi(\varphi_g^*(x))$. 
Here, $\psi$ denotes the absolute value function and $\partial_{\varphi}\psi(e)$ is the subdifferential of $\psi$ with respect to the variable $\varphi$ evaluated at $e\in \Re$. See \eqref{eq:subdifferential-absolute-value} for an expression of the subdifferential $\partial_{\varphi}\psi(\cdot)$ from which we know that $s$ satisfies $|s|\leq 1$. 

%%%%%%%%
\noindent Manipulating the equations above yields 
$$ 
\begin{array}{ll}
	 z_g^*(x) &= x-\gamma s a \\
	 \varphi_g^*(x) &= a^\top z_g^*
	(x)+b-\dfrac{\gamma}{\lambda}s\\
	& =a^\top x+b-\gamma\left(\frac{1}{\lambda}+\left\|a\right\|_2^2\right)s.
\end{array}
$$
It follows that 
\begin{equation}\label{eq:opt-zg-xig}
\left\{	\begin{aligned}
		& z_g^*(x) = x-\gamma s a \\
		& \varphi_g^*(x) =\eta \left(\rho(x)-s\right),
	\end{aligned}\right.
\end{equation}
with the notations $\rho(x)$ and $\eta$ defined in \eqref{eq:rho(x)}. 
To determine the value of $s$ (which is, as defined before, a function of $\varphi_g^*(x)$), we now consider three cases: 
\begin{itemize}
	\item If $\rho(x)>1$, then $\varphi_g^*(x)>0$, $s=1$ and 
$$ 
\left\{\begin{aligned}
	& z_g^*(x) = x-\gamma a \\
	& \varphi_g^*(x) =\eta \left(\rho(x)-1\right)
\end{aligned}\right.
$$
\item If $\rho(x)<-1$, then $\varphi_g^*(x)<0$, $s=-1$ and 
$$ 
\left\{\begin{aligned}
	& z_g^*(x) = x+\gamma a \\
	& \varphi_g^*(x) =\eta \left(\rho(x)+1\right)
\end{aligned}\right.
$$
\item If $\left|\rho(x)\right|\leq 1$, then we have necessarily $\varphi_g^*(x)=0$ (assuming the contrary leads to a contradiction) and hence $s=\rho(x)$ in this case. As a consequence,  we obtain from \eqref{eq:opt-zg-xig} that 
$$ 
\left\{\begin{aligned}
	& z_g^*(x) = x-\gamma \rho(x) a \\
	& \varphi_g^*(x) =0
\end{aligned}\right.
$$
\end{itemize}
Combining these three cases give the expressions in \eqref{eq:closed-form-zg}.  
%$s=1$ if $\rho(x)>1$, $s=-1$ if $\rho(x)<-1$, $s=\rho(x)$ if $\left|\rho(x)\right|\leq 1$
\end{pf}
%%%%%

\begin{rem}
Note that the function $z_g^*$ in \eqref{eq:(z,xi)} can be equivalently expressed in term of the proximal operator as 
\begin{equation}
	z_g^*(x) =  \prox_{g_1}(x)
\end{equation}
with $g_1$ being defined on $\Re^n$ by 
$$
\begin{aligned}
		g_1(z)& =\inf_{\varphi\in \Re}g(z,\varphi)\\ %= \Sat_{\gamma/\lambda}\left(a^\top z+b\right).  
		  & = \frac{\gamma^2}{\lambda}\left[\dfrac{1}{2}\Sat_1^2\left(\dfrac{a^\top z+b}{\gamma/\lambda}\right)+\left|\dz_1\left(\dfrac{a^\top z+b}{\gamma/\lambda}\right)\right|\right]
\end{aligned}
$$
with $\Sat_1^2(\cdot)$ being understood as the square of the saturation function. 
\end{rem}
%%%%%%%%%%
By using Lemma \ref{lem:lasso-type}, it is immediate to see that the pair $(z_t^{(i)},\varphi_{ti}^{(i)})$ in \eqref{eq:relaxation-Lasso}  can be expressed as 
\begin{align}
&z_t^{(i)} = z_t^{(i-1)}+\gamma_{ti} \Sat_{1} (\rho_{it})W_t^2c_{ti} \label{eq:z-lasso}\\ 
&\varphi_{ti}^{(i)}  = \eta_{it}\dz_1(\rho_{it}), \label{eq:xi-lasso}
\end{align}
where $\rho_{it}\!= \!\frac{1}{\eta_{it}}(y_{ti}-c_{ti}^\top z_t^{(i-1)})$ and $\eta_{it}\!= \!\gamma_{ti}(\frac{1}{\lambda_{ti}}+\left\|W_tc_{ti}\right\|_2^2)$.

\subsection{Some other robustness-inducing loss functions}
We note here that many other loss functions can be used in \eqref{eq:measurement-update} to design a robust estimator. 
Examples of such functions include the Huber, Vapnik and Log-abs  loss functions (see \cite{Barron2019,Sayin14-TSP} for the latter). We derive for these three functions closed-form update laws when they are used in the component-wise relaxation as in \eqref{eq:closed-form-Abs} (See  Lemmas \ref{lem:Huber}--\ref{lem:Vapnik-Loss} in the appendices). The  proximal observers associated with all these  loss functions have  been evaluated in simulation, see Section \ref{sec:simulation}. It should be noted that Vapnik's loss function is not positive definite  as demanded by Assumption \ref{A:positive-definite}. Therefore it may not achieve convergence of the error to zero even when the measurement are noise-free.

\section{Illustrative simulations}\label{sec:simulation}
\subsection{A linear example}
To illustrate the performance of the estimation paradigm discussed in the previous sections, we consider an example of LTI system with one input, $2$ outputs and $3$ state variables in the form \eqref{eq:system} with 
$f_t(x)=Ax+Bu_t$ and $C_t=C$ where
\begin{equation}\label{eq:linear-example}
	\begin{aligned}
		&A=\begin{pmatrix}-1  &   1  &   0 \\
			-1  &   0 &    0 \\
			0  &  -1  &  -1 \\
			\end{pmatrix}, \quad 
		B = \begin{pmatrix}
		-1\\ 0\\ 0
		\end{pmatrix}, \quad C = \begin{pmatrix}
		1&0&0\\ 0&0&1
		\end{pmatrix}. %\\ 
		%%
		%&C = \begin{pmatrix}
		%1&0&0\\ 0&0&1
		%\end{pmatrix}
	\end{aligned}
\end{equation}
The initial state is set to $x_0=(\begin{matrix} 10& 5& 5\end{matrix})$ and the input $u_t$ is selected to be a sine function in the form $u_t=\sin(2\pi\nu_0 tT_s)$ where $\nu_0=0.1$Hz, $T_s=0.1$. 
Note that here, the matrix $A$ has all its eigenvalues on the unit circle which implies that $A$ is not Schur stable. As a consequence, it is not possible to construct a stable open-loop observer (i.e., a simple simulator of the dynamic equation independently of the measurement) for this example. Therefore, the design of a valid state estimator requires the use of measurements (even though they are potentially erroneous).  For this purpose we note that $(A,C)$ is observable.   

\begin{figure}[h!]%
\centering
\subfloat[]{\includegraphics[scale=.77]{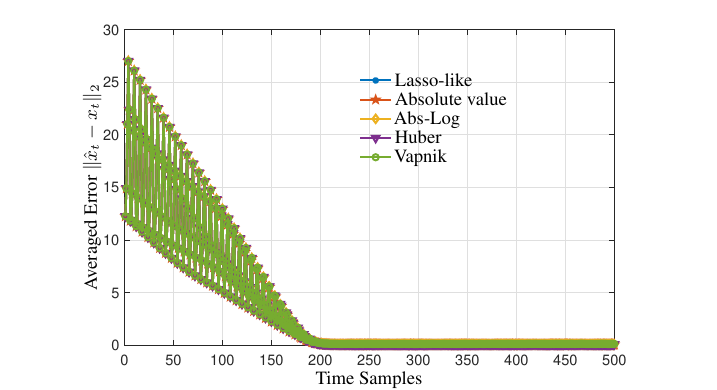}} \\
\subfloat[a Zoom on the time interval $\interval{450}{500}$.]{\includegraphics[scale=.77]{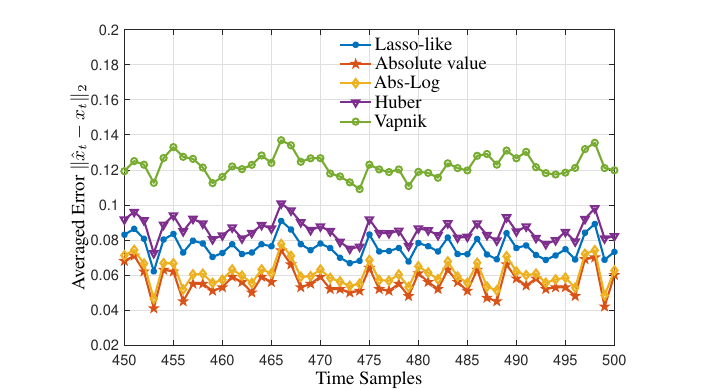}}%
\caption{Averaged norms (over $100$  realizations) of the absolute estimation errors when only the sparse measurement is active. Minimum dwell-time between consecutive occurrences of nonzero instances of the sparse noise is fixed $5$ samples (component-wise).  }%
\label{fig:SPARSE:Averaged-Error-versus-Time}%
\end{figure}

\begin{figure}[h]%
\centering
\includegraphics[scale=.8]{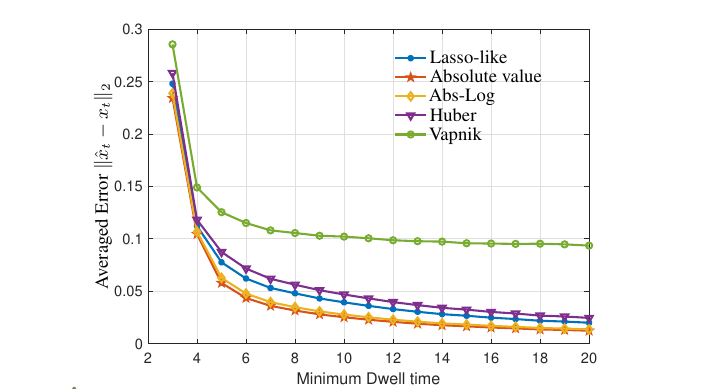}%
\caption{Averaged norms,  over $100$  realizations, of the (absolute) estimation errors when only the sparse measurement is active. }%
\label{fig:Error-State-DwellTime}%
\end{figure}

\paragraph{Effect of the impulsive noise}
In a first experiment, we set  all dense noises to zero, i.e.,  $w_t=0$ and $\nu_t=0$ for all $t$ and consider only the effect of the impulsive noise $\left\{\zeta_t\right\}$. It is sampled from a white Gaussian noise with distribution $\mathcal{N}(0,100I_2)$. To make this sequence sparse we randomly set to zero some of its entries $\zeta_{ti}$ over time in such a way that two consecutive nonzero instances of $\left\{\zeta_{ti}\right\}_{i,t}$ are separated by a minimum dwell time of $\tau_{\dwell}=5$ samples. 
We implement all the 5 adaptive estimators\footnote{Only the component-wise relaxation versions are implemented.} discussed in the paper with the parameters recorded in Table \ref{tab:parameters}. 
\begin{table}%
\centering
\begin{tabular}{|l|l|}
\hline
Estimator & Parameters \\ \hline
Lasso & $\lambda=2$, $\gamma=0.1$  \\ \hline
Absolute value &  $\lambda=0.1$ \\ \hline
Abs-Log  & $\lambda=0.1$, $\mu=1000$ \\ \hline
Huber  & $\lambda=0.1$, $\mu=0.08$ \\ \hline
Vapnik & $\lambda=0.1$, $\epsilon=0.07$\\ \hline
\end{tabular}
\caption{Parameters of the robust state estimators as implemented. For all estimators, the weighting matrix taken to be $W_t=I$. }
\label{tab:parameters}
\end{table}  
It turns out from the norm of the state estimation error displayed in Figure \ref{fig:SPARSE:Averaged-Error-versus-Time} that all algorithms perform very well in attenuating the effect of the impulsive noise on the estimate. The errors converge to  reasonably small values despite a high frequency of nonzero occurrences in $\left\{\zeta_{ti}\right\}$. The price to pay however for this robustness property is a slow convergence rate. This is because the parameters (in particular the threshold levels of the component-wise relaxed versions) of the estimator are kept small and fixed at all times. Selecting those parameters dynamically would accelerate the convergence.      
Figure \ref{fig:Error-State-DwellTime} shows the average error (after convergence) for different values of the the dwell time $\tau_{\dwell}$. Naturally, the estimation error gets lower as the dwell time increases. 

To further robustify the estimators against impulsive noise, one can consider a two-stage estimation process as follows: \\
(a) Based on the estimate $\hat{\zeta}_{ti}$ (of the $i$-th component of the impulsive noise $\zeta_t$) provided by the robust estimators at each time, remove the "bad" measurements whenever the estimate $\hat{\zeta}_{ti}$ exceeds some appropriately chosen thresholds $T_{ti}$ defined for example, by 
$$T_{ti}=\max\big(\min\big(T_{t-1,i},\left|y_{ti}-c_{ti}^\top\hat{x}_t\right|\big),\epsilon_0\big), \: \:  \epsilon_0=0.01$$ 
with $\hat{x}_t$ being the estimate returned by the robust proximal observer.  \\
(b) apply a standard observer, for example, a steady-state Kalman filter, on the good measurements. \\
Implementing such an estimation scheme  does not seem to affect the convergence properties in this example. Empirical evidence shows that it considerably reduces the steady-state estimation error.  In the case of the Absolute Value and the Abs-Log proximal estimators, the steady state estimation error even reaches (numerical) zero despite the presence of the sparse noise.  

%%%%%%
%%%%%%%
%%%%%%%%%%%%%%%%%%%

\begin{figure}%
\centering
\includegraphics[scale=.77]{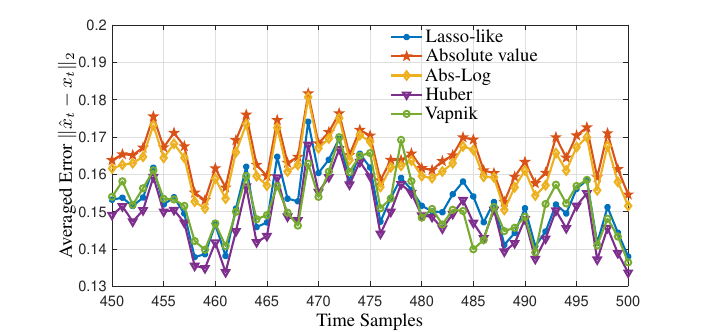}%%width=.55\textwidth,height=5cm
\caption{Averaged norms  of the absolute estimation errors (over $100$  realizations) when all process and measurement noises are active. A zoom on the interval $\interval{450}{500}$.}%
\label{fig:SPARSE-DENSE:Averaged-Error-versus-Time}%
\end{figure}

\begin{figure}%
\centering
\includegraphics[scale=.7]{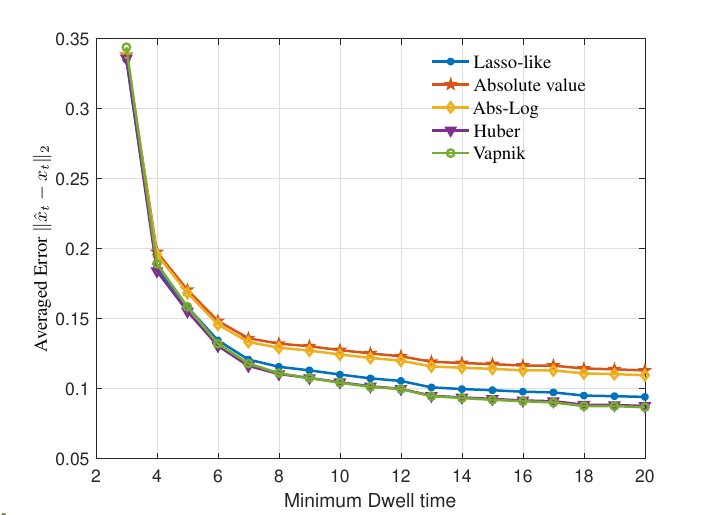}%
\caption{Averaged norms of the (absolute) estimation errors when all process and measurement noise are active over $100$  realizations. }%
\label{fig:AvgError-versus_DwellTime}%
\end{figure}

\paragraph{Combined effect of both impulsive and dense noises}
To illustrate robustness to both dense and sparse noises, we carry out a second experiment where  the dense noises $\left\{w_t\right\}$ and $\left\{\nu_t\right\}$ are now independently sampled from uniform distributions with support $\interval{-0.1}{0.1}$  while the impulsive noise $\left\{\zeta_t\right\}$ is still active. 
Assuming the same estimators' parameters as in the first experiment (see Table \ref{tab:parameters}), we obtain the results displayed in Figure \ref{fig:SPARSE-DENSE:Averaged-Error-versus-Time} in term of the absolute estimation error on the interval $\interval{450}{500}$ (for a better visualization purpose). As expected there is a degradation of the estimation errors as a consequence of the action of the dense noises. Also the bad data detection threshold appears to be less effective than in the first case.  But the estimation errors remain reasonably bounded despite the large values taken by the impulsive noise and the relatively high frequency at which its nonzero values occur over time (See Figure \ref{fig:AvgError-versus_DwellTime}).  
%%%
%%%%%%%%%%%%%%%%%%%%%%%%%%%
%%%%

\subsection{A nonlinear example}
\begin{figure}[h!]
\centering
\includegraphics[scale=.75]{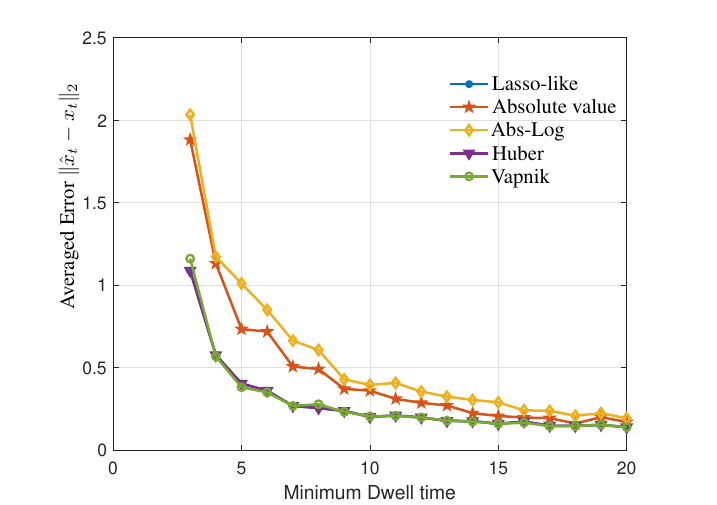}%
\caption{Nonlinear example  -- Averaged norms of the (absolute) estimation errors when all process and measurement noise are active over $100$  realizations.}%
\label{fig:AvgError-versus_DwellTime-NL}%
\end{figure}

We consider now a nonlinear system of the form \eqref{eq:system} 
with $f_t(x)=Ax+Bu_t+F(x)$ and $y_t=c^\top x_t$, where $A$ and $B$ are defined as in \eqref{eq:linear-example}, $c^\top=[1\:  1 \:  1]$ and $F$ is a nonlinear vector field given by  
$
F(x)=\big[\begin{matrix}
	 \sin(x_1+x_2), & 
	 \sin x_1 \cos x_2, &
	 \Sat_1(x_3)
\end{matrix}\big]^\top. 
$
Considering a similar experiment as in the first example, we implement the detection-correction strategy. Contrarily to the linear case where a Kalman steady state filter was used in the correction step, here, the same estimator is used for both detection and correction steps.  
The resulting absolute estimation errors are depicted in Figure \ref{fig:AvgError-versus_DwellTime-NL}. What this shows is that for  all the implemented five instances of our estimation framework, the errors follow the same trend as in the case of the  linear example. 
%%%%

\section{Conclusion}
In this paper we have discussed a framework for designing state estimators which are robust to impulsive measurement noise.  
At each discrete time step, the observer performs a prediction and a measurement update, the latter being formulated as an optimization problem.  We showed that robustness properties are induced  by an appropriate choice of the objective function which is minimized for the measurement update. In particular, we considered five particular choices of loss functions and derive, under some relaxations, analytic update  expressions of the associated estimators. Future work will concern a more complete performance analysis of the observer. One open question of interest is how to systematically and optimally select online the parameters of the algorithms in order to achieve a better performance. A desirable property would be one guaranteeing that the estimation error is stabilized while remaining completely insensitive to the extreme values of the noise.

%%%%%%%%%%%%%%%%%%%%%%%%%%%%%%%%%%%%%%%%%%%%%%%%%%%%%%%%%%
\appendix

\section{Uniform complete observability (UCO)}\label{subsec:UCO}
This appendix complements the arguments of Section \ref{subsec:kalman} concerning the Kalman filter. We prove that if the open loop system $\left\{A_k,C_k\right\}$ is UCO, then so is the closed-loop $\left\{\bar{A}_k,C_k\right\}$ with  $\bar{A}_k=M_{k+1}^{-1}A_k$. 
For this purpose, note that $M_{k+1}^{-1}=I_n-L_{k+1}C_{k+1}$,  where 
$L_k=P_{k|k-1}C_k^\top\left(R_k+C_kP_{k|k-1}C_k^\top\right)^{-1}C_k.$
By letting $F_k=I_{n_y}-C_kL_k$, it is not hard to show that $F_k= R_k\left(R_k+C_kP_{k|k-1}C_k^\top\right)^{-1}$ and hence that the sequence 
$\left\{F_k\right\}$ is uniformly bounded under the assumptions of Section \ref{subsec:kalman}. Lemma \ref{lem:equivalence-UCO} can therefore be applied to obtain the claim. 
\begin{lem}\label{lem:equivalence-UCO}
Let $\left\{(A_k,C_k)\right\}$ and  $\{(\bar{A}_k,C_k)\}$ be two sequences in $\Re^{n\times n}\times \Re^{n_y\times n}$ satisfying $\bar{A}_k=(I_n-L_{k+1}C_{k+1})A_k$ for some sequence $\left\{L_k\right\}$ in $\Re^{n\times n_y}$ such that $\{F_k\triangleq I_{n_y}-C_kL_k\}$ is uniformly bounded  in the sense that $\underline{\sigma} I_{n_y}\preceq F_k^\top F_k\preceq\overline{\sigma}I_{n_y}$ for some strictly positive $\underline{\sigma},\overline{\sigma}$. Then  $\left\{(A_k,C_k)\right\}$ is UCO if and only if  $\left\{(\bar{A}_k,C_k)\right\}$ is UCO. 
\end{lem}
%%%%%
This lemma was indeed stated in \cite[Lemma A.2]{Moore80-IJC} but only a concise sketch of the proof  was given.  We provide here a more detailed proof for completeness. 
\begin{pf}
Let $m>0$ be an integer. For $t\geq m$, denote with $\mathcal{N}_{t,t-m}$ and $\overline{\mathcal{N}}_{t,t-m}$  the observability grammians of $\left\{(A_k,C_k)\right\}$ and  $\left\{(\bar{A}_k,C_k)\right\}$ respectively on the finite  time horizon $[t-m,\: t]$. 
Let $\overline{\Gamma}_{t,t-m}$ and $\Gamma_{t,t-m}$ be the $(m+1)$-horizon observability matrices associated to both pair sequences. For illustration, $\overline{\Gamma}_{t,t-m}$ has the form     
\begin{equation}\label{eq:obsv-matrix}
\overline{\Gamma}_{t,t-m}=\begin{pmatrix}C_{k}\\C_{k+1}\bar{A}_{k}\\ C_{k+2}\bar{A}_{k+1}\bar{A}_{k}\\ \vdots \\ C_t\bar{A}_{t-1}\cdots \bar{A}_{k}\end{pmatrix}
\end{equation}
where $k=t-m$ for compactness. 
Then direct calculations show that  $\overline{\Gamma}_{t,t-m} = S_{t,t-m}\Gamma_{t,t-m}$ where  $ S_{t,t-m}$ is  a block-lower triangular matrix with a structure given by
$$
S_{t,t-m} \!= \!\!
\left(\!\begin{array}{llllll}
I_{n_y}\\
0 & F_{k+1}\\
0 & S_{t,k}^{32}& F_{k+2}\\
0  & S_{t,k}^{42} & S_{t,k}^{43} & F_{k+3}\\
\vdots & \cdots & \cdots &  \cdots & \ddots\\
0 &  \cdots & \cdots  & \cdots  & \cdots & F_{t} 
\end{array}\!\right)
$$
where the off-diagonal blocks take the following forms  
$$
\begin{aligned}
	&S_{t,k}^{32} = -F_{k+2}C_{k+2}A_{k+1}L_{k+1}\\
	&S_{t,k}^{42} =   -F_{k+3}C_{k+3}A_{k+2}\bar{A}_{k+1}L_{k+1} \\
	&S_{t,k}^{43} = -F_{k+3} C_{k+3}A_{k+2}L_{k+2}. 
\end{aligned}
$$
Using this remark we can write 
\begin{equation}\label{eq:equality-grammians}
	\begin{aligned}
		\overline{\mathcal{N}}_{t,t-m}& =\overline{\Gamma}_{t,t-m}^\top \overline{\Gamma}_{t,t-m}\\
		& = {\Gamma}_{t,t-m}^\top (S_{t,t-m}^\top  S_{t,t-m}){\Gamma}_{t,t-m}.
	\end{aligned}
\end{equation}
From the assumption of the lemma, we know that $S_{t,t-m}$ is bounded. 
Recall from Weyl's inequalities \cite[Chap. 5]{Bernstein09-Book}  that   
$$ |\det(S_{t,t-m})| = \prod_{i=1}^{q}|\lambda_i(S_{t,t-m})|=\prod_{i=1}^{q}\sigma_i(S_{t,t-m}).$$
with $q=n_y(m+1)$ and $\lambda_i$ and $\sigma_i$ referring to eigenvalues and singular values respectively with the indexation being such that $|\lambda_1(S_{t,t-m})|\geq|\lambda_2(S_{t,t-m})|\geq\ldots>|\lambda_q(S_{t,t-m})|$ and $\sigma_1(S_{t,t-m})\geq \sigma_2(S_{t,t-m})\geq\ldots \geq \sigma_q(S_{t,t-m})$. 
Now assume for contradiction that  $\left\{S_{t,t-m}\right\}$ is not uniformly bounded away from $0$, i.e.,  $\inf_{t}\sigma_q(S_{t,t-m})=0$. 
Using the boundedness of $\left\{S_{t,t-m}\right\}$, this implies that  
$\inf_{t}|\det(S_{t,t-m})|=0$ which in turn implies that $\inf_{t}|\det(F_t)|=0$ and $\inf_{t}\sigma_{n_y}(F_t)=0$, which constitutes a  contradiction with the assumption of uniform boundedness of $\left\{F_t\right\}$. We therefore conclude that there necessarily exists a constant $s_0>0$ such that $(S_{t,t-m}^\top  S_{t,t-m})\succeq s_0 I_q$. Returning to \eqref{eq:equality-grammians} we see that the claim of the lemma holds.  
\end{pf}

\subsection{Some technical results}
In this appendix we study three other possible candidate loss functions to be used in \eqref{eq:gt} for achieving robust measurement update in the state estimation framework discussed in this paper.
\begin{table*}[h!]
\centering
\small
 \renewcommand{\arraystretch}{2.2}
\begin{tabular}{|l|l|l|l|}
\hline
Loss & $g(z)$     & $\prox_{g}(x)$ & Ref. \\ \hline 
Absolute value & $\alpha|a^\top z+b|$ &  $ x-\alpha\Sat_{1}\left(\dfrac{a^\top x+b}{\alpha \left\|a\right\|_2^2}\right) a $ & Lemma \ref{lem:absolute-value}\\ \hline
Lasso-type & $\begin{array}{l}\inf_{\varphi\in \Re}g(z,\varphi)\\ \mbox{(see Eq. \eqref{eq:g(z,xi)})}\end{array}$ &  $x-\gamma\Sat_{1}\left(\dfrac{a^\top x+b}{\gamma(1/\lambda+\left\|a\right\|_2^2)}\right) a$ & Lemma \ref{lem:lasso-type} \\ \hline
Huber & $\alpha h_{\mu}(a^\top z+b)$ &  $ x-\alpha\Sat_1\left(\dfrac{a^\top x+b}{\mu+\alpha \left\|a\right\|_2^2}\right) a$ & Lemma \ref{lem:Huber} \\ \hline 
Log-Abs & $\alpha L_{\mu}(a^\top z+b)$ & $x-\dfrac{\alpha\mu \omega}{1+s\mu \omega} a$ & Lemma \ref{lem:Logorithmic-loss}\\ \hline
Vapnik & $\alpha \max(|a^\top z+b|-\epsilon,0)$  & $x-\alpha  d_{\sigma,\epsilon}(a^\top x+b) a$ & Lemma \ref{lem:Vapnik-Loss}\\ \hline
\end{tabular}
\caption{Some robust proximal operators.}
\end{table*}

\begin{lem}[Update with Huber loss]\label{lem:Huber}
Let $g:\Re^n\rightarrow\Re$ be the function defined by 
$g(z) = \alpha h_{\mu}(a^\top z+b) $ 
with  $a\in \Re^n$,  $b\in \Re$, $\alpha\in \Re_+$ and $h_\mu:\Re\rightarrow\Re_+$ being the Huber loss defined, for  a given $\mu>0$, by  
\begin{equation}\label{eq:Huber-function}
	h_\mu(e)=\left\{
	\begin{array}{ll}
	\dfrac{1}{2\mu}e^2 & \mbox{if }|e|\leq \mu\\
	|e|-\dfrac{\mu}{2} & \mbox{otherwise.}
	\end{array}
	\right.
\end{equation}
Then 
\begin{equation}
	\prox_{g}(x) = x-\alpha\Sat_1\left(\dfrac{a^\top x+b}{\mu+\alpha \left\|a\right\|_2^2}\right) a,
\end{equation}
with $\Sat_1$ defined as in \eqref{eq:Saturation}. 
\end{lem}
%%%%
\begin{pf}
$h_\mu$ is convex and so for $g$ which is therefore subdifferentiable. Consequently, $\prox_{g}$ is well-defined. 
We note that 
$$z^*=\prox_{g}(x)=\argmin_{z}\left[\dfrac{1}{2}\left\|z-x\right\|_2^2+ g(z)\right]$$
 if and only if $0\in z^*-x+ \partial g(z^*)=z^*-x+\alpha a r^*$ with $r^*\in \partial h_\mu(a^\top z^*+b)$. 
Because $\partial h_\mu(e)=\big\{\Sat_1(\frac{e}{\mu})\big\}$, we obtain 
$$r^*= \Sat_1\left(\dfrac{a^\top x+b-\alpha r^*\left\|a\right\|_2^2}{\mu}\right).$$
Solving this equation for $r^*$ yields
$$r^*= \Sat_1\left(\dfrac{a^\top x+b}{\mu+\alpha \left\|a\right\|_2^2}\right).$$
Finally, inserting this in the expression of $z^*$ gives the claimed result.  
\end{pf}
%%%%%%%%%%%%%%%%%%%

\begin{lem}[Update with Abs-log loss]\label{lem:Logorithmic-loss}
Let $g:\Re^n\rightarrow\Re$ be the function defined by  $g(z)=\alpha L_{\mu}(a^\top z+b)$ 
with $a\in \Re^n$, $b\in \Re$, $\alpha\in \Re_+$ and  $L_\mu:\Re\rightarrow\Re_+$ defined by 
\begin{equation}\label{eq:Logarithmic-function}
	L_{\mu}(e)=|e|-\dfrac{1}{\mu}\ln(1+\mu|e|)
\end{equation}
for a given $\mu>0$.  
Then 
\begin{equation}
	\prox_g(x)=x-\dfrac{\alpha\mu \omega}{1+s^\star\mu \omega} a,
\end{equation}
with $s^\star=\sign(a^\top x+b)\in \left\{-1,0,+1\right\}$ and 
\begin{equation}
	\begin{aligned}
		&\omega=\dfrac{r+s^\star\sqrt{\Delta}}{2\mu}\\
		& r=\mu(a^\top x+b)-s^\star(1+\alpha\mu \left\|a\right\|_2^2)\\
		& \Delta = r^2+4\mu |a^\top x+b|.
	\end{aligned}
\end{equation}
\end{lem}
%%%%%%%%%%%
%%%%%%%
\begin{pf}
Note, to begin with, that the function $L_\mu$ in \eqref{eq:Logarithmic-function} is convex. To see this,  observe that $L_\mu(e)$ can be expressed as $L_\mu(e)=h(|e|)$ where $h$ is given by $h(e)=e-\frac{1}{\mu}\ln(1+\mu e)$. Since $h$ is convex on $\Re_+$ (which can be checked by noting that the second derivative of $h$ is strictly positive) and nondecreasing, $L_\mu$ is convex by the composition rule \cite[Chap. 3]{Boyd04-Book}.  
It follows that $g$ as defined in the statement of the lemma is also convex and hence subdifferentiable. \\ 
Now $z^*=\prox_{g}(x)$ if and only if 
$0\in z^*-x+\partial g(z^*)=z^*-x+\alpha a \partial L_\mu(a^\top z^*+b)$. 
After observing that 
$$\partial L_\mu(e)=\frac{\mu |e|}{1+\mu |e|}s(e),$$
with $s$ being a set-valued map defined by $s(e)=\partial|e|$ (see \eqref{eq:subdifferential-absolute-value}), we get  
$$0\in z^*-x+\alpha a \frac{\mu |\omega^*|}{1+\mu |\omega^*|} s(\omega^*),  $$
where we have set $\omega^* = a^\top z^*+b$.  
We distinguish the following three cases according to the sign of $\omega^*$: 
\begin{itemize}
	\item Case 1: $\omega^* >0$. Then $s(\omega^*)=\left\{1\right\}$ and hence  
$$
\begin{aligned}
	&\left\{
	\begin{aligned}
		&\omega^* >0\\
		&0 = z^*-x+\frac{\alpha\mu  \omega^*}{1+\mu \omega^*}a
	\end{aligned}
		\right.\\
\end{aligned}
$$
$$
\begin{aligned}
	 \Leftrightarrow & 
		\left\{
	\begin{aligned}
		&\omega^* = a^\top x+b-\frac{\alpha\mu  \omega^*}{1+\mu \omega^*}\left\|a\right\|_2^2>0\\
		& z^*=x-\frac{\alpha\mu  \omega^*}{1+\mu \omega^*}a
	\end{aligned}
	\right.
\end{aligned}
$$

\item Case 2: $\omega^* <0$. Then $s(\omega^*)=\left\{-1\right\}$ and hence   
$$
\begin{aligned}
	&\left\{
	\begin{aligned}
		&\omega^* <0\\
		&0 = z^*-x+\frac{\alpha\mu  \omega^*}{1-\mu \omega^*}a
	\end{aligned}
		\right.
		\\
		\Leftrightarrow \quad 
		&	\left\{
	\begin{aligned}
		&\omega^* = a^\top x+b-\frac{\alpha\mu  \omega^*}{1-\mu \omega^*}\left\|a\right\|_2^2<0\\
		& z^*=x-\frac{\alpha\mu  \omega^*}{1-\mu \omega^*}a
	\end{aligned}
	\right.
\end{aligned}
$$

\item Case 3:  $\omega^* =0$. Then $s(\omega^*)=\interval{-1}{1}$ and hence  
$$
\left\{
\begin{aligned}
	&\omega^* =0\\
	&0 = z^*-x
\end{aligned}
	\right.
	\quad \Leftrightarrow \quad 
		\left\{
\begin{aligned}
	&\omega^* = a^\top x+b=0\\
	& z^*=x
\end{aligned}
	\right.
$$
\end{itemize}
Solving for $\omega^*$ in these cases while taking into account the associated constraints and the fact that a solution necessarily exists, yields the result. 
\end{pf}
%%%%%%%%%%%%%%%
\begin{lem}[Update with Vapnik loss]\label{lem:Vapnik-Loss}
Let $g:\Re^n\rightarrow\Re_+$ be the function defined by $g(z)=\alpha \psi_\epsilon(a^\top z+b)$, with $\psi_{\epsilon}(e)=\max(|e|-\epsilon,0)$ for some $\epsilon\geq 0$, $a\in \Re^n\setminus\left\{0\right\}$, $b\in \Re$, $\alpha>0$. 
Then
$$\prox_g(x) = x-\alpha a d_{\sigma,\epsilon}(a^\top x+b), $$
with $\sigma=\epsilon+\alpha \left\|a\right\|_2^2$
and 
%%%%%
\begin{equation}\label{eq:dsigma}
	d_{\sigma,\epsilon}(e) = 
	\left\{\begin{array}{ll}
	\sign(e), & |e|> \sigma\\
	0, & |e|< \epsilon\\
	\dfrac{1}{\sigma-\epsilon}\big(e-\epsilon\sign(e)\big), & \epsilon\leq \left|e\right|\leq \sigma.
	\end{array}\right. 
\end{equation}
\end{lem}
%%%%%%%%%%%%%%%%%%
\begin{pf}
\noindent $z^*=\prox_{g}(x)$ if and only if 
$0\in z^*-x+\partial g(z^*)=z^*-x+\alpha a \partial \psi_\epsilon(a^\top z^*+b)$. 
By recalling that
$$
\partial \psi_\epsilon(e)
=
\left\{\begin{array}{ll}
\left\{\sign(e) \right\}& |e|> \epsilon \\
\left\{0\right\} & |e|<\epsilon\\
\interval{0}{1} & e=\epsilon\\
\interval{-1}{0} & e=-\epsilon,\\
\end{array}\right. 
$$
we get that $z^*$ satisfies 
$$
z^* = 
\left\{\begin{array}{lll}
x-\alpha a \sign(a^\top z^*+b), &  |a^\top z^*+b|> \epsilon\\
x, & |a^\top z^*+b|< \epsilon\\
x-\alpha a \lambda_+, &a^\top z^*+b= \epsilon\\
x-\alpha a \lambda_-, &a^\top z^*+b= -\epsilon\\
\end{array}\right. 
$$
with $(\lambda_-,\lambda_+) \in \interval{-1}{0}\times \interval{0}{1}$. \\
In the first case, if we let $s^*= \sign(a^\top z^*+b)$, then we get 
$$\begin{aligned}
	|a^\top z^*+b|&=\left|\phi(x)-\alpha\left\|a\right\|_2^2s^*\right|\\
	&=s^*\left(\phi(x)-\alpha\left\|a\right\|_2^2s^*\right)>\epsilon,
\end{aligned} $$
with $\phi(x)=a^\top x+b$. 
This is equivalent to $s^*\phi(x)> \epsilon+\alpha\left\|a\right\|_2^2$, which in turn translates into $s^*=\sign(\phi(x))$ and  $\left|\phi(x)\right|> \epsilon+\alpha\left\|a\right\|_2^2$.  It follows that $z^*=x-\alpha a \sign(\phi(x))$ for $|\phi(x)|> \epsilon+\alpha\left\|a\right\|^2=\sigma$. \\
The second case gives immediately $z^*= x$ for $|\phi(x)|< \epsilon$. \\
The third and fourth case give
$$
\left\{\begin{array}{ll}
\lambda_+= \dfrac{\phi(x)-\epsilon}{\alpha \left\|a\right\|^2 }, &  \epsilon\leq  \phi(x)\leq \sigma\\
\lambda_- = \dfrac{\phi(x)+\epsilon}{\alpha \left\|a\right\|^2 }, &  -\sigma\leq  \phi(x)\leq -\epsilon. 
\end{array}\right.
$$
To recapitulate,  $z^*=\prox_{g}(x)$ is given by 
$$
z^* = 
\left\{\begin{array}{ll}
x-\alpha a \sign(\phi(x)), & |\phi(x)|> \sigma\\
x, & |\phi(x)|< \epsilon\\
x-\dfrac{a}{\left\|a\right\|^2}\big(\phi(x)-\epsilon\sign(\phi(x))\big), & \epsilon\leq \left|\phi(x)\right|\leq \sigma.
\end{array}\right. 
$$
This is the result stated by the lemma. 
\end{pf}

%%%%%%%%%%%%%%%%%%%%%%%%

\bibliographystyle{abbrv}

\end{document}